\documentstyle[fullpage, 11pt]{article}

\def\SBIMSMark#1#2#3{
 \font\SBF=cmss10 at 10 true pt
 \font\SBI=cmssi10 at 10 true pt
 \setbox0=\hbox{\SBF Stony Brook IMS Preprint \##1}
 \setbox2=\hbox to \wd0{\hfil \SBI #2}
 \setbox4=\hbox to \wd0{\hfil \SBI #3}
 \setbox6=\hbox to \wd0{\hss
             \vbox{\hsize=\wd0 \parskip=0pt \baselineskip=10 true pt
                   \copy0 \break%
                   \copy2 \break%
                   \copy4 \break}}
 \dimen0=\ht6   \advance\dimen0 by \vsize \advance\dimen0 by 8 true pt
                \advance\dimen0 by -\pagetotal
 \dimen2=\hsize \advance\dimen2 by .25 true in
     \setbox0=\hbox to 3.1 true in{
                \vbox to \ht6{\hsize=3 true in \parskip=0pt  \noindent  
  {\it  Invent. Math.}~{\bf 122} (1995), 35--69
                \vfill}}
  \ht0=0pt \dp0=0pt
 \ht6=0pt \dp6=0pt
 \setbox8=\vbox to \dimen0{\vfill \hbox to \dimen2{\copy0 \hss \copy6}}
 \ht8=0pt \dp8=0pt \wd8=0pt
 \copy8
 \message{*** Stony Brook IMS Preprint #1, #2 ***}
}
\newcommand{\1}{{{\mathchoice {\rm 1\mskip-4mu l} {\rm 1\mskip-4mu l}
{\rm 1\mskip-4.5mu l} {\rm 1\mskip-5mu l}}}}

\newcommand{\R}{{\bf R}}

\newcommand{\Z}{{\bf Z}}
\newcommand{\C}{{\bf C}}

\newcommand{\p}{{\partial}}
\newcommand{\al}{{\alpha}}

\newcommand{\be}{{\beta}}

\newcommand{\Om}{{\Omega}}
\newcommand{\om}{{\omega}}
\newcommand{\eps}{{\varepsilon}}
\newcommand{\de}{{\delta}}

\newcommand{\ga}{{\gamma}}
\newcommand{\Ga}{{\Gamma}}
\newcommand{\ka}{{\kappa}}
\newcommand{\la}{{\lambda}}
\newcommand{\si}{{\sigma}}

\newcommand{\T}{{\bf T}}

\newcommand{\Cc}{{\cal C}}

\newcommand{\Ll}{{\cal L}}

\newcommand{\Nn}{{\cal N}}

\newcommand{\Pp}{{\cal P}}

\newcommand{\Uu}{{\cal U}}

\newcommand{\Si}{{\Sigma}}
\newcommand{\Ham}{{\rm Ham }}
\newcommand{\IMP}{\Longrightarrow}
\newcommand{\can}{{\rm can}}
\newcommand{\Totvar}{{\rm Totvar}}
\newcommand{\area}{{\rm area}}
\newcommand{\vol}{{\rm vol}}
\newcommand{\id}{{\rm id }}
\newcommand{\im}{{\rm Im }}

\newcommand{\MS}{{\medskip}}

\newcommand{\NI}{{\noindent}}
\newcommand{\proof}[1]{\noindent{\bf Proof#1:\  }}
\newcommand{\jdef}[1]{{\bf #1}}
\newcommand{\QED}{\hfill$\Box$\medskip}
\newcommand{\Cal}{{\rm Cal\,}}

\newcommand{\Tilde}{\widetilde}
\newtheorem{theorem}{Theorem}[section]
\newtheorem{cor}[theorem]{Corollary}

\newtheorem{definition}[theorem]{Definition}
\newtheorem{remark}[theorem]{Remark}
\newtheorem{lemma}[theorem]{Lemma}

\newtheorem{prop}[theorem]{Proposition}
\newtheorem{proposition}[theorem]{Proposition}
\newcommand{\at}{{@}}

\title{Hofer's $L^{\infty}$\/-geometry:\\ energy and stability  of Hamiltonian
flows,\\ part II}
\author{Fran\c{c}ois Lalonde\thanks{Partially supported by NSERC grant
OGP 0092913
and FCAR grant ER-1199.} \\ Universit\'e du Qu\'ebec \`a Montr\'eal
\\ (flalonde\at math.uqam.ca) \and Dusa McDuff\thanks{Partially supported by
NSF grant DMS 9103033 and NSF Visiting Professorship for Women GER
9350075.} \\ State University of New York at Stony Brook \\
(dusa\at math.sunysb.edu)}
\date{}

\begin{document}
\maketitle

\SBIMSMark{1995/3a}{February 1995}{}
\thispagestyle{empty}

\MS

\section*{Abstract}

In this paper we first show that the necessary condition introduced in our
previous paper
is also a sufficient condition for a
path to be a geodesic in the group $\Ham^c(M)$ of  compactly supported
Hamiltonian symplectomorphisms. This applies with no restriction on $M$.
We then discuss conditions which guarantee that such a path
minimizes the Hofer length.  Our argument relies on a general geometric
construction (the gluing of monodromies) and on an extension of
Gromov's   non-squeezing theorem both to more general manifolds
and to more general capacities.  The manifolds we consider
are quasi-cylinders, that is spaces
homeomorphic to $M \times D^2$ which are symplectically ruled over $D^2$. 
When we work with the usual capacity (derived from embedded balls),
we can prove the existence of paths which
minimize the length among all homotopic paths, provided that
$M$ is semi-monotone.  (This restriction occurs because of the
well-known difficulty with the theory of $J$-holomorphic curves in
arbitrary $M$.)  However, we can only prove the
existence of length-minimizing paths (i.e. paths which minimize length
amongst {\it all} paths, not only the homotopic ones) under
even more
restrictive conditions on $M$, for example when $M$ is exact and convex
or  of dimension $2$.    The new 
difficulty is caused by the possibility that there are non-trivial and
very short loops in $\Ham^c(M)$.  When such length-minimizing paths do exist,
we can extend the Bialy--Polterovich calculation of the Hofer norm on a
neighbourhood of the  identity ($C^1$-flatness). 

Although it applies to a more
restricted class of manifolds,  the Hofer-Zehnder 
capacity seems to be better adapted to
the problem at hand, giving sharper estimates in many situations.
Also the capacity-area inequality for split cylinders extends more easily to 
quasi-cylinders in this case. As applications,
we generalise Hofer's estimate of the time for which
an autonomous flow is length-minimizing to some manifolds other
than $\R^{2n}$, and derive new
results such as the unboundedness of Hofer's metric on some closed
manifolds, and a linear rigidity result.

\section{Statement of main results}

\subsection{Geodesics}

In this paper, unless specific mention is made, $M$ will be any symplectic
manifold, compact or not, with or without boundary, but always geometrically
bounded at infinity in the non-compact case. When $\p M \neq \emptyset$,
all Hamiltonians have compact support in $M - \p M$. 
 
For the convenience of the reader, we begin by recalling some of the basic
definitions from Part I (\cite{LALM1}).  The length $\Ll(\ga)$ of the path
$\ga = \phi_{t\in[a,b]}$ in the group $\Ham^c(M)$ generated by the 
compactly supported Hamiltonian
function $H_t$ is defined by:
$$
\Ll(\ga) = \Ll(\phi_t) = \int_a^b \|H_t\|dt,
$$
where
$$
\|H_t\| = \Totvar(H_t) = \sup_{x\in M} H_t(x) - \inf_{x\in M}
H_t(x).
$$
The Hofer norm $\|\phi\|$ of $\phi\in \Ham^c(M)$ is defined to be the infimum of
the lengths of all paths from the identity to $\phi$.

We say that the path $\ga = \phi_{t\in [a,b]}$ is
\jdef{regular} if its tangent vector $\dot \phi_t$ is non-zero  for all $t\in
[a,b]$. A \jdef{stable geodesic} $\ga = \phi_{t\in [a,b]}$ is a 
regular  path
which is a local minimum for $\Ll$ on the space $\Pp$ of all
paths  $\psi_{t\in[a,b]}$ homotopic to $\phi_{t\in [a,b]}$
with fixed endpoints $\phi_a,\phi_b$, where $\Pp$ is given the usual
$C^{\infty}$ topology.

As in Part I, we will say that a path $\ga = \phi_{t\in[0,1]} $ has
a certain property {\bf at each moment} if every $s\in [0,1]$ has a closed
connected neighbourhood $\Nn(s)$ in $[0,1]$ such that, for all
subintervals $\Nn$ containing $s$, the subpath $\phi_{t\in \Nn}$ has this
property.  Very often (but not always)  it will suffice to check that the
property in question holds on the interval $\Nn(s)$ itself since it will then
automatically hold on any subinterval.  For example, this is clearly the case
for the property of being a stable geodesic.  
  We recall that a regular path $\ga$ is said to be a
\jdef{geodesic} if it is a stable geodesic at each moment.  Hence each $s\in
[a,b]$ has a connected neighbourhood  $\Nn(s)$ such that $\phi_{t\in \Nn(s)}$
is a stable geodesic.
 Unless specific mention is made to the contrary, we will
normalise $t$ so that it takes values in the interval $[a,b]=[0,1]$, and will
assume that $\phi_0 = \1$.   Also, often it is convenient to talk about a path 
of
Hamiltonians $H_{t\in[0,1]}$ rather than the corresponding path in
$\Ham^c(M)$.

A point $P$ such that $H_t(P) = \sup_{x\in M} H_t(x)$ for all $t\in [0,1]$ is
called a {\bf fixed maximum} of the corresponding  path $\ga$ (or
equivalently of $H_{t\in [0,1]}$.)  A similar definition applies to a fixed
minimum $p$. We often denote these by $P,p$, and write $q$ for a fixed
extremum.   Borrowing terminology from \cite{BP}, we will often call a
path with a fixed maximum and
minimum {\bf quasi-autonomous}.  Clearly, every autonomous 
(i.e. time-independent) Hamiltonian $H$ is quasi-autonomous.

In Part I we established the following necessary condition for a regular path
to be a geodesic.

\begin{theorem} A regular path $\phi_t, t\in I$, in $\Ham^c(M)$ is a geodesic
only if its generating Hamiltonian has at least one fixed
maximum and one fixed minimum at each moment. 
\end{theorem}

We will prove below in \S3.1 that:

\begin{theorem}\label{thm:geodsuff} For any manifold $M$, the above necessary
condition is also sufficient. \end{theorem}

When $M = \R^{2n}$ the above result can be deduced from the work of
Bialy--Polterovich in~\cite{BP}.  We refer the reader to the other papers of this series \cite{LALM1,LALM3}
for the characterization of the stability of geodesics.

\subsection{Length-minimizing properties of geodesics}

In order to state our results on the existence of length-minimizing paths,
we need to introduce the idea of the {\bf capacity} $c(H)$  of a
Hamiltonian $H_{t\in [0,1]}$.  We will explain this precisely in \S2 below.
However, roughly speaking, it is, for each given choice $c$ of a capacity, the
minimum of the $c$\/-capacity of the regions over and under the graph of 
$H_{t\in [0,1]}$. 

    We recall that, for any symplectic manifold $M$, with or without boundary,
compact or not, the Gromov capacity of $M$ is 
$$
c_G(M) = \sup \{a : \mbox{there is a symplectic embedding} \; f:B^{2n}(a) \to 
(M - \p M) \}
$$
where $B^{2n}(a)$ denotes the standard closed ball of $\R^{2n}$ of capacity 
$a = \pi r^2$. The Hofer-Zehnder capacity $c_{HZ}(M) \in [0, \infty]$ of 
$M$ is the supremum
over all positive real numbers $m \in [0, \infty)$ such that there is a
surjective function $H:M \to [0,m]$ equal to $m$ on $\p M$ and outside some 
compact set of $M$, and whose Hamiltonian flow has no non-constant closed
trajectory in time less than $1$.

Thus, for instance, if $c = c_G$ is the Gromov capacity, $c_G(H)$ is 
the maximum capacity of a symplectic ball
which embeds on both sides of the graph of $H$.   One can
show that if $H_t$ is sufficiently $C^2$-small, then $c(H) \geq \Ll(H_t)$
for any choice of $c$.

We also need to put a
condition on $M$ in order to use the theory of $J$-holomorphic curves.
Following the terminology of~\cite{DUSAL}, we will say that $(M,\om)$ is  {\bf
weakly monotone}  if for every spherical homology class $B\in H_2(M)$
$$
\om(B) > 0,\;\;c_1(B) \ge 3-n \quad\IMP \quad c_1(B)\ge 0.
$$
This condition is satisfied if
either $\dim M \le 6$ or $M$ is {\bf  semi-monotone}, i.e. there  is a constant
$\mu \ge 0$ such that, for all spherical homology classes $B\in H_2(M)$,
$$
c_1(B) = \mu \om(B).
$$
As explained in~\cite{DUSAL}, weak monotonicity is exactly the condition
under which the theory of $J$-holomorphic curves behaves well.  Since we
apply this theory to the product $M\times S^2$ rather than to $M$ itself, our
arguments do not work for all weakly monotone $M$.   However they do
work when $\dim M \le 4$, since this implies that $M\times S^2$ is weakly
monotone and also, by a slightly different argument, when $M$ is
semi-monotone.  Thus it applies to the semi-monotone manifolds $M = \T ^{2n}$
(the standard torus)
and $M = \C P^n$ (complex projective space) as well as to some but not all
products of the form $\C
P^k\times \C P^\ell$.  

Recall that a manifold $M$ is called weakly exact 
if $[\om]\mid_{\pi_2(M)} = 0$.

\subsection*{Paths minimal in a homotopy class}

\begin{theorem}\label{thm:monotone} 
\begin{description}
\item[(i)]
 Let $M$ be have dimension $\le 4$ or be
semi-monotone. Let $H_{t \in [0,1]}$ be any path with
$$
c_G(H) \geq \Ll(H_t).
$$
Then $H_{t \in [0,1]}$ is length-minimizing amongst all paths homotopic
rel endpoints to $H_{t \in [0,1]}$. \\
\item[(ii)] The same conclusion holds if $c_{HZ}(H) \geq \Ll(H_t)$ and
$M$ is weakly exact.
\end{description}
\end{theorem}

\begin{remark}\rm (i) $\;$ As a consequence of the Non-Squeezing Theorem, $\Ll(H_t)$
is actually the maximal possible value of $c_G(H)$. This also true of $c_{HZ}(H)$
when $M$ is weakly exact as a consequence of the $c_{HZ}$-area inequality, and of
any capacity when $M = \R^{2n}$. See \S2.
Thus the theorem states in fact that the path is length-minimizing in its
homotopy class as soon as $c(H)$ reaches its maximal value.

\NI
(ii) $\;$ Note that in this theorem, we do not need any other
hypothesis on $H_{t \in [0,1]}$. In particular, this theorem implies
that any path with $c_G(H) = \Ll(H_t)$ must be a stable geodesic and
therefore must have two fixed points by the necessary condition for stability.
A more elaborate argument leads to linear rigidity, see Theorem~\ref{thm:linear-rig}
at the end of \S5.
\end{remark}

This theorem is proved in \S5.  (A simplified proof is available when $\dim\, M =
2$: see Remark~\ref{rmk:split}.)  To apply it we need to know when
$c(H) \geq \Ll(H_t)$.  This seems a hard problem in general.  However, 
the inequality $c(H) \geq \Ll(H_t)$ is always true when $H$ is $C^2$-small and,
as we shall see in \S3, the inequality $c_{HZ}(H) \geq \Ll(H_t)$ is also true when
$H$ is autonomous and has no non-constant closed trajectory
in time $< 1$. The same statement holds for $c_G$ when $M$ has dimension $2$.   
Hence we find:

\begin{cor} \label{cor:HS-homotopic} Let $M$ be a weakly exact manifold or any
surface. 
Let $H_{t \in [0,1]}$ be an autonomous Hamitonian whose flow has no
non-constant closed trajectory in time less than $1$. Then the flow $\phi_{t \in [0,1]}$
generated by $H$
is length-minimizing among all homotopic paths rel endpoints.
\end{cor} 

\subsection*{Short loops and length-minimizing paths}

Our methods do not allow us to compare  the
lengths of non homotopic  paths in $\Ham^c(M)$ for arbitrary $M$.  However,
there is a topological reason for this. Consider the length function on the
fundamental group of $\Ham^c(M)$: $$
\Ll:\pi_1(\Ham^c(M))\to [0,\infty), \quad
\Ll([\ga]) = \inf_{\ga\in[\ga]} \Ll(\ga).
$$
We will say that $\Ham(M)$ has {\bf short loops} if $\{0\}$ is not
an isolated point in the image of $\Ll$.  Thus, if $M$ has short loops,
there is
for each $\eps > 0$  a class $[\ga]$ in $\pi_1(\Ham^c(M))$ which can be
represented by a loop of length $< \eps$ but not by a loop of length $<
\eps/2$.  When $M$ does not have short loops we define  $r_1 = r_1(M) > 0$ to be
$$
r_1 = \inf\, \left( \{\im\, \Ll:\pi_1 \to \R\} \cap (0, \infty) \right)
$$
when this is not empty, and $\infty$ otherwise.  If $M$ has short loops
we put $r_1 = 0$.

As an example,  observe that $\Ham^c(\R^{2n})$ does
not have short loops, since any given loop may be homotoped to one which
is arbitrarily short by conjugating it by an appropriate rescaling factor.
A similar argument shows:

\begin{lemma}\label{le:shloop} Suppose that $(M,\om)$ is exact and is convex 
in the sense that it admits a  contracting Liouville vector field
whose flow exists for all positive time.
Then $\Ham^c(M)$ does not have short loops.
\end{lemma}

Another situation in which $M$ obviously cannot have short loops is when
$\pi_1(\Ham^c(M))$ is finite or cyclic.  This applies to surfaces: it is well-known
that when $M$ has dimension $2$, $\Ham^c(M)$ is contractible, except in the
case of the $2$-sphere where 
$$ 
\pi_1\Ham(S^2) = \pi_1(SO(3)) = \Z/2. 
$$
It also applies to certain $4$-manifolds, for example to $\C P^2$ with its
standard form and to $S^2\times S^2$, provided that the latter manifold has a
product form $\si_0\oplus\si_1$ which has the same integral on both sphere
factors.  This follows from Gromov's calculation of the homotopy type of
$\Ham^c(M)$ in these cases: see~\cite{G}.

When $(M,\om)$ does not have short loops, it is sometimes possible to calculate
the invariant $r_1$.  For example, 
when $M = \C P^2$ or is any surface other than the sphere, $\pi_1(\Ham^c(M)) =
0$ and so  $r_1 = \infty$.  Less trivially, in \S5.2 we  use
Theorem~\ref{thm:monotone} to show that:

\begin{lemma}\label{le:turn}  Let $\phi_{t \in [0,1]}$ be an essential
loop in $\Ham(S^2)$, where $S^2$ has area $A$.
Then $\Ll(\phi_{t \in [0,1]}) \ge A$.
\end{lemma}

Observe that the path generated by a multiple of the height
function, which rotates the sphere about the north-south axis through a
complete turn, is an essential loop with length exactly $A$.  Hence

\begin{cor} If $(S^2, \om)$ has area $A$, then $r_1(S^2) = A$.
\end{cor}
\MS

The following theorem collects together our main results on the existence 
of length-minimizing paths.

\begin{theorem}\label{thm:lengmin}  Let  $M$ have
dimension $\le 4$ or be semi-monotone, and have no short loops.  In the
non-compact case, assume also that it has bounded geometry at infinity. 
Then the path
$H_{t \in [0,1]}$ is length-minimizing amongst all paths with the same endpoints
in the following three cases:
\begin{description}
\item[(i)]  when $c_G(H) = \Ll(H_t) \le r_1 /2$; 

\item[(ii)]  when $c_{HZ}(H) = \Ll(H_t) \leq r_1 /2$ and $M$ is weakly exact;

\item[(iii)] when $M = \R^{2n}$ and there is some capacity for which
$$
c(H) \geq \Ll(H_t).
$$ 
\end{description}
\end{theorem}

\proof{}  By Theorem~\ref{thm:monotone}, any path $\ga$ with
$c_G(H) = \Ll(H_t)$ is
is length-minimizing among all homotopic paths.  Suppose
its length is $\le r_1/2$, and let $\ga'$ be some other path
with the same endpoints and with $\Ll(\ga') < \Ll(\ga)$.
Then the loop $-\ga' * \ga$ has length $< r_1$ and so represents a
loop which has an arbitrarily short representative.  We may choose
a representative $\la$ of length $\eps$ where
$$
\Ll(\ga') + \eps < \Ll(\ga).
$$
The $\ga'*\la$ is homotopic to $\ga$ and is shorter: a contradiction.
This proves (i).  Essentially the same proof applies in the other two cases.
\QED

\begin{cor}\label{cor:HS-absolute} Let $M$ be an exact and convex manifold
($\R^{2n}$ for instance) or any surface,
not necessarily compact, and let $H$ be
an autonomous Hamiltonian with time-$1$ map $\phi$. Suppose that
$H$ has no non-constant closed trajectory in time less than $1$, and if
$M$ is a sphere, that $\|H\| \le \frac 12 \area(M)$. Then
$$
\|\phi\| = \|H\|.
$$
\end{cor}

    This generalises Hofer's criterion in \cite{HOF}. See Richard~\cite{RICH}
for further results in this direction. It is possible that our
methods can be used in the non-autonomous case too: then the
above theorem would be true for any quasi-autonomous Hamiltonian
satisfying the above hypothesis.

\subsection*{Diameter of $\Ham^c(M)$}

 Much of our work in
studying $\Ham^c(M)$ for a closed orientable surface was prompted by the 
question of
whether this group has bounded diameter.  If $M$ is open and of finite volume,
 the
existence of the Calabi invariant implies that $\Ham^c(M)$ has infinite 
diameter.
However, in view of Sikorav's result in~\cite{S} 
that the subgroup $\Ham^c(B^{2n}(1))$ has
bounded diameter in $\Ham^c(\R^{2n})$, it is not clear what the answer should
be for compact $M$.  We have not been able to decide this when  $M = S^2$.
When the genus is $> 0$, one can use the energy-capacity inequality to
show:

\begin{prop}  Any closed orientable surface of genus $> 0$ supports an
autonomous Hamiltonian $H$  whose flow $\{\phi_t^H\}_{t\ge 0}$ is length
minimizing for all $t$. \end{prop}

\begin{cor}\label{cor:unbdd} For these manifolds $M$, the group  $\Ham(M)$ has
unbounded diameter with respect to the Hofer norm.
\end{cor}

This corollary can be extended to all manifolds of the form $M \times \Sigma$,
where $\Sigma$ is any closed orientable surface of genus at least $1$ and
$M$ any manifold. This uses holomorphic methods and will be developed 
elsewhere.

\subsection*{Local flatness}

If $H_{t\in [0,1]}$ has both a fixed maximum and minimum and is also
sufficiently small, it follows from a version of the Hamilton-Jacobi
equation that $\Ll(H_t)$ may be measured by the difference between the
actions of its fixed maximum and minimum.  We therefore can extend the
local flatness result of Bialy--Polterovich~\cite{BP} in this case:

\begin{cor} \label{cor:BP} If $(M,\om)$ is semi-monotone and has no
short loops, there is a $C^1$-neighbourhood
$\Uu$ of the identity in $\Ham^c(M)$ which, when given the Hofer metric,  is
isometric to a neighbourhood of $\{0\}$ in a Banach space.  
\end{cor}

\subsection*{Properties holding at each moment}

The above results concern the length of whole paths.  One can also
look at properties of pieces of paths.  By 
definition, a geodesic is stable at each moment, and it is natural to 
wonder how much this property can be strengthened.  Ideally, one would like
geodesics to be length-minimizing at each moment.  However, we have not been
able to prove this for arbitrary $M$.  The next result is an easy consequence of
Theorem~\ref{thm:monotone} and the equality $c_G(H) = \Ll(H_t)$ for
$C^2$-small Hamiltonians.

\begin{prop}\label{prop:moment1}  Suppose that $(M,\om)$ is
semi-monotone.  Then every geodesic
minimizes length among homotopic paths  at each
moment.  Moreover, if in addition $(M,\om)$ has no short loops, every geodesic
is length-minimizing at each moment.  \end{prop}

This means that for such $M$ our
geodesics behave like geodesics in Riemannian geometry.  In particular, any
quasi-autonomous Hamiltonian isotopy is absolutely length-minimizing for some
interval of time $[0,t_0], t_0 > 0$, a result proved by Bialy--Polterovich
in~\cite{BP} when $M = \R^{2n}$ using other methods.
\MS

The final result we will mention here concerns rational 
manifolds $(M,\om)$, that is to say those which
have a positive {\bf  index of rationality} $r$.  Here $r$ is defined to be the
largest positive element of the extended reals $(0,\infty]$
  such that $\om(\pi_2(M))\subset r\Z$, if this exists, and to be
$0$ otherwise.

\begin{prop}\label{prop:moment2}  Suppose that $(M,\om)$ has bounded
geometry at infinity, and that it has index of rationality $r>0$.  Then, at
each moment, a geodesic minimizes length among paths which are
homotopic to it by paths of length $<2r$.  More precisely, for each $s\in [0,1]$
there is an interval $\Nn(s)$ such that, whenever $s\in \Nn\subset \Nn(s)$,
the path $\phi_{t\in \Nn}$ has minimal length among all paths which are
smoothly homotopic to it through paths of length $< 2r$. 
\end{prop}

A version of the above Proposition which applies to smooth deformations  of
certain long initial paths is given in Theorem~\ref{thm:rational}.

\subsection{Non-squeezing results for quasi-cylinders}

In order to establish the length-minimizing properties of
geodesics we need an extension of the non-squeezing theorem to \lq\lq
quasi-cylinders".  We will state our basic results on quasi-cylinders in this
subsection, and show in \S2 how they are related to Hofer's geometry.

    Given a capacity $c$, a symplectic manifold $(M, \om)$, and $a \in 
(0, \infty)$, we say that the cylinder $M \times D^2(a)$ endowed with the
split form, satisfies the $c$\/-area inequality if 
$$
c(M \times D(a)) \leq a.
$$
Note that, by the definition of a capacity, the space $\R^{2n} \times D^2(a)$
satisfies the $c$\/-area inequality for any $c$ and any $a > 0$ (see 
\cite{LAL} for instance for the definition of a (intrinsic) capacity).

We now state the capacity-area inequalities. The first one, due to
Gromov and Lalonde-McDuff, applies to all manifolds. The second one,
due to Floer, Hofer and Viterbo, applies to rational manifolds in a
specific range. 

\begin{theorem}[Non-Squeezing Theorem,~\cite{G,LALMCD}]\label{thm:ns}
The split cylinder $M \times D^2(a)$ satisfies the $c_G$-area inequality
for any $(M, \om)$ and any $a > 0$.
\end{theorem}

As explained in~\cite{LALMCD}, there is a direct proof using
$J$-holomorphic curves when $M$ is \lq\lq nice" (eg if
$M$ has dimension $\le 4$ or is semi-monotone), and a much more elaborate
indirect proof for arbitrary $M$.

\begin{theorem}[HZ-capacity-area inequality,~\cite{FHV,HV}] \label{thm:HZ-ineq}
Let $M$ be a rational manifold of index of rationality $r \in (0, \infty]$.
Then $M \times D^2(a)$ satisfies the $c_{HZ}$-area inequality whenever 
$a \leq r$.
\end{theorem}

\begin{definition}\label{qc} \rm Let $(M,\om)$ be a symplectic manifold
and $D$ a set diffeomorphic to a disc in $(\R^2, \si = du \wedge dv)$.
Then the manifold $Q = (M\times D, \Om)$ endowed with a symplectic
form $\Om$ is called a \jdef{quasi-cylinder} if
\begin{description}
\item[(i)] $\Om$ restricts to $\om$ on each fiber $M\times pt$;

\item[(ii)] $\Om$ is the
product $\om \oplus  \si$ near the boundary $M \times \p D$,
and, in the case when $M$ is non-compact, outside a set of the form
$X\times D$, for some compact subset $X$ in $ M$.
\end{description}

 A quasi-cylinder $Q$ is said to be \jdef{split} if there is a
symplectomorphism
$$
(M\times D, \Om) \;\to \;(M\times D, \om \oplus \si)
$$
which is the identity near the boundary.  Its \jdef{area} is the number
$A$ such
that
$$
\vol\,(X\times D, \Om) = A \, \vol\,(X,\om),
$$
where $X$ is a compact set as in (ii) above.
\end{definition}

We will see in Lemma~\ref{le:area} below that the area of $Q$ does not depend on
the choice of $X$.  It is not known
whether every quasi-cylinder is split, though by the results
of~\cite{RULED} on the structure of ruled surfaces  this is the case when
$\dim M = 2$.

We will see in \S2 that every pair of paths from $\1$ to $\phi$ gives rise to
 two quasi-cylinders, whose areas are related to the lengths of the paths.
When the paths are sufficiently $C^2$-close, these quasi-cylinders
split and we can use the usual capacity-area inequalities to obtain results valid
for arbitrary $(M,\om)$.  This is the method used to prove
Theorem~\ref{thm:geodsuff}, for example.   In general, we must deal with
arbitrary quasi-cylinders, and our results will hold only for certain $(M,\om)$.

For ease of exposition, we now state  slightly simplified
versions of the results in \S4.  We will say that the $c$\/-area inequality
holds for a quasi-cylinder $Q$ of area $A$ if $c(Q) \leq A$ (when $c = c_G$,
we will also call this a non-squeezing theorem).

\begin{prop}\label{prop:quasi-ns0} Let $M$ be a symplectic manifold with
bounded geometry and $Q = (M \times D, \Om)$ be a
quasi-cylinder of area $A$.
\begin{description}
\item[(i)] Then the non-squeezing theorem holds for $Q$
whenever $M$ has dimension $\le 4$ or is semi-monotone. \\
\item[(ii)] The $c_{HZ}$-area inequality holds for $Q$ whenever $M$
is weakly exact.
\end{description}
\end{prop}

The proof of (i) is a somewhat tricky extension of the usual proof of the
non-squeezing theorem (as given  for example in~\cite{LALMCD}), and
involves a non-embedding result for  a pair consisting of a ball and a
cylinder. But, not surprisingly, the proof of (ii) is much easier: indeed
the relation between $c_G$ and $c_{HZ}$ is somewhat similar to the relation 
between handle decompositions and Morse functions. The $c_{HZ}$-capacity
is more manageable since functions are easier to homotope than their level sets.
  Since this proposition is the fundamental result which underlies
the proof of Theorem~\ref{thm:monotone} it holds for the same $M$ as does
that theorem:  see Remark~\ref{rmk:condit}.

The following result on manifolds with positive index
of rationality is a deformation result analogous
to~\cite[Lemma~3.3]{LALMCD}.

\begin{prop}\label{prop:deformation0}
Let $\Om_{s \in [0,1]})$ be a smooth
$1$\/-parameter family of symplectic forms on $M\times D$ such that
each manifold  $Q_s = (M \times D, \Om_s)$ is a quasi-cylinder. Assume
that $Q_0$ is split.
\begin{description}
\item[(i)]   If $\area \, (Q_s) < 2r$ for all $s$,
then the non-squeezing theorem holds for $Q_1$;

\item[(ii)]
 If there is a smooth symplectic isotopy
$$
g_s: B^{2n+2}(c) \to Q_s, \quad 0 \le s \le 1
$$
which at time $s=0$ has the form
$$
B^{2n+2}(c) \hookrightarrow B^{2n}(c)\times D(c)\stackrel{f\times g}
{\longrightarrow} M\times D
$$
where $c \ge 2r$, and if
$$
\area \,Q_s < c + 2r, \quad 0 \le s \le 1,
$$
then $\area \, Q_1 \ge c$.
\end{description}
\end{prop}

\subsection{Organization of the paper}

This paper is organised as follows.  In \S2 we explain the connection between
Hofer's geometry and quasi-cylinders and work out the simplest properties of
quasi-cylinders.  In \S3 we describe some methods of embedding balls, and apply
them to prove Theorem~\ref{thm:geodsuff} and results in dimension $2$ such as
Lemma~\ref{le:turn}. We also construct Hamiltonians which realise the 
Hofer-Zehnder capacity. \S4 is the heart of the paper, and concerns non-squeezing
theorems for quasi-cylinders.  All applications to Hofer's geometry 
are proved in \S5.

We wish to acknowledge the hospitality of the Newton Institute 
where this paper was finalized.

\section{Hofer's geometry and quasi-cylinders}

In this section we first explain how the lengths of paths are related to
the capacities of quasi-cylinders.  In \S2.2 we develop the basic
properties of quasi-cylinders.

\subsection{Gluing monodromies}\label{ss*:gluing}

    The basic geometric construction needed 
consists in defining a symplectic space (a ``quasi-cylinder'') 
made from two Hamiltonian
isotopies which are homotopic rel endpoints. If $H_{t \in [0,1]},
K_{t \in [0,1]}$ are the corresponding Hamiltonians, we attach the lower
part of the graph of one of these maps to the upper part of the graph of the 
other by the  identification of the graphs which match their 
characteristic foliations.

Up until now, we have insisted that  geodesics $\ga$ be regular paths, which 
means that
their time-derivatives $\frac {\p}{\p t}\ga(t) $ vanish nowhere. In this 
section it will be convenient for technical reasons to consider paths which are
constant  for
$t$ in some intervals $[0,\eps]$, $[1-\eps, 1]$ near the endpoints, and  have
non-vanishing time-derivative elsewhere.  Note that any regular path can
be reparametrized so that this is the case, without changing its length. 
Further, since we will also only consider paths which have at least one fixed
maximum and minimum, we may  normalise their generating
Hamiltonians $H(x,t)$ by assuming that  $\inf_x H(x,t) = 0$ for all
$t$.\footnote[1]{Note that in general the function $t\mapsto \inf_x H(x,t)$, 
though
continuous, need not be smooth.} 

We will write $\Gamma_H$ for the 
graph of $H$: 
$$
\Gamma_H = \{(x, H(x,t), t)\in M \times \R \times [0,1]\}.
$$
Observe that, if $M\times \R\times [0,1]$ is provided with
the symplectic form $\Om = \om \oplus ds\wedge dt$, then the
lines of the characteristic flow induced on the hypersurface $\Ga_H$ are
precisely
$$
(\phi_t(x_0), H(\phi_t(x_0)), t), \quad 0\le t \le 1.
$$
We write $R_H^-, R_H^+$ respectively for the region under and over the graph:
\begin{eqnarray*}
R_H^- & = & \{(x, s,t) \in M \times \R\times [0,1]: 0 \leq s \leq H(x,t)\} \\
R_H^+ & = & \{(x, s,t) \in M \times \R\times [0,1]: H(x,t) \leq s \leq 
\max_x H(x,t)\}
\end{eqnarray*}
For a small parameter $\nu > 0$, we denote by $R_H^-(\nu), R_H^+(\nu)$
the thickened regions:
\begin{eqnarray*}
R_H^-(\nu/2) &=& \{(x, s,t) \in M \times \R\times [0,1]: \la(t) \leq s 
\leq H(x,t)\} \\
R_H^+(\nu/2) &=& \{(x, s,t) \in M \times \R\times [0,1]: H(x,t) \leq s \leq 
\mu_H(t)\}
\end{eqnarray*}
where $\la(t)$ is a smooth function taking values in $[-\de, 0]$
which equals $-\de$ except near $0,1$, and where 
$\mu_H(t)$ is a smooth function such that $m_H(t) =\max_x H(x,t)
\le \mu_H(t) \le \max_x H(x,t) + \de$.  We assume that $\la(t)$ and
$\mu_H(t)$ are chosen so that
$$
U_H = \{(s,t): \la(t) \le s\le \mu_H(t)\} \subset \R^2
$$
has a smooth boundary. Then $\nu$ refers to the area added in the thickening:
$$
\nu/2 = \int_0^1 - \la(t) \, dt \;=\;\int_0^1 (\mu_H(t) - 
\max_x H(x,t)) \, dt.
$$
Of course,
$$
R_H(\nu) = R_H^-(\nu/2) \cup R_H^+(\nu/2) = M \times U_H
$$
is the direct product of $M$ with a disc of area $\|H\| + \nu$.
Note that the two halves $R_H^{\pm}(\nu/2)$ of
$R_H(\nu)$ have corners at $t = 0,1$: their boundaries are the union of
two smooth pieces joined along $M\times \{0\}$ and $M\times \{1\}$.
One of these boundary  pieces is a product and has trivial holonomy,
and the other is $\Ga_H$ with holonomy $\phi_1$.

Suppose that $\{\phi_t\}$ and $\{\psi_t\}$ are two isotopies with
$\phi_1 = \psi_1$, and let $H_t, K_t$ be the corresponding Hamiltonians.
Then we may join $R_K^+(\nu/2)$ to $R_H^-(\nu/2)$
by the map $\Ga_K \to \Ga_H$ given by:
$$
(x, s, t)  \mapsto (\phi_t\circ \psi_t^{-1}(x),
         s - K(x) + H(\phi_t\circ \psi_t^{-1}(x)),\, t).
$$
Observe that this map is a symplectomorphism on the whole of
$R_K^+(\nu/2)$.  Indeed, if $f_t$ is the isotopy generated by $F_t$, then
$$
(y,s,t)\mapsto (f_t(y), s + F_t(f_t(y)), t)
$$
is a symplectomorphism.  Here, we have
$$
f_t = \phi_t
\circ\psi_t^{-1},\quad F_t(y) = H_t(y) -
K_t(\psi_t\circ\phi_t^{-1}(y)). 
 $$
Thus
$$
R_{H,K}(\nu) =  R_H^-(\nu/2) \cup R_K^+(\nu/2) = 
\{(x, s,t): \la(t) \leq s \leq  \mu_K(t)
 + F_t(x) \}
$$
is a manifold with smooth boundary.

Because $\phi_1 = \psi_1$, the
monodromy round $\p R_{H,K}(\nu)$ is trivial.   Therefore, there is a
symplectomorphism 
$$
\Phi_K:  M \times {\rm nbhd}( \p U_H) \to {\rm nbhd}( \p R_{H,K}(\nu))
$$
given by
$$
(x,u,v) \mapsto (\phi_v\circ\psi_v^{-1}(x), u - \mu_H(v) + 
F_v(\phi_v\circ\psi_v^{-1}x), v)
$$
near the part of the boundary where $u = \mu_H(v)$, and by the identity
elsewhere.  In fact, this formula makes sense on the whole of the thickening
$R_{H,K}(\nu) - R_{H,K}$ of $\p R_{H,K}$.   Because the loop
$\{\phi_t\circ\psi_t^{-1}\}$ contracts in $\Ham(M)$, $\Phi_K$ extends to a
diffeomorphism 
$$
\Phi_K: M\times U_H \to R_{H,K}(\nu)
$$
of the form
$$
(x,u,v)\mapsto (f_{u,v}(x),\,  s(x,u,v), \,v),
$$
where each $f_{u,v}\in \Ham^c(M)$.  It follows easily  that
$(M\times U, \Phi_K^*(\Om))$ or, equivalently $(R_{H,K}(\nu),\Om)$, is a
quasi-cylinder.

The following very simple lemma is the key to our results.

\begin{lemma}\label{le:compare}  Suppose that $\Ll(K_t) < \Ll(H_t) = A$.
Then, for sufficiently small $\nu > 0$, at least one of the quasi-cylinders
$(R_{H,K}(\nu),\Om)$  and $(R_{K,H}(\nu),\Om)$ has area $< A $.
\end{lemma}

\proof{} Choose $v>0$ so that 
$$
\Ll(K_t) + 2\nu < \Ll(H_t),
$$
and suppose first that $M$ is compact.
Evidently,
\begin{eqnarray*}
\vol\,(R_{H,K}(\nu)) + \vol\,(R_{K,H}(\nu)) & = &\vol\,(R_H(\nu)) + 
\vol\,(R_{K}(\nu)) \\
& = & (\vol\,M).(\Ll(H_t) + \Ll(K_t) + 2\nu)\\
&< & 2(\vol\,M).\Ll(H_t).
\end{eqnarray*}
If $M$ is non-compact, we may restrict to a large compact piece  $X$  of $M$
and then take the volume.\QED

In order to use this proposition, we consider the following. Let $c$ be any capacity
and define the capacity $c(H)$ of 
the Hamiltonian $H_t$ 
as the minimum of $\inf_{\nu > 0} c(R_H^-(\nu/2))$ and $\inf_{\nu > 0} c(R_H^+(\nu/2))$.
\footnote{This means that, for instance,
the $c_{HZ}$-capacity of
$H_{t \in [0,1]}$ is defined by looking at autonomous Hamiltonians defined
under (or over) the graph of a non-autonomous Hamiltonian.
}  
Suppose now that 
$$
c(H) \geq \Ll(H_t).
$$
(Actually, $c_G(H)$ cannot be greater than $\Ll(H_t)$ by
Theorem~\ref{thm:ns}, since $R_H^-(\nu/2)\cup R_H^+(\nu/2)$ is a split
cylinder of area $\Ll(H_t)+\nu$ for all $\nu$. The same remark applies for $c_{HZ}$ when
$M$ is weakly exact by Theorem~\ref{thm:HZ-ineq}). Then for any Hamiltonian $K_{t \in
[0,1]}$ which generates an isotopy homotopic rel endpoints to $\phi_t$, there
are symplectic embeddings $$
   R_H^-(\nu/2) \hookrightarrow R_{H,K}(\nu) \quad {\rm and} \quad
   R_H^+(\nu/2) \hookrightarrow R_{K,H}(\nu).
$$
Therefore $c(R_{H,K}(\nu)) \geq \Ll(H_t)$ and
$c(R_{K,H}(\nu)) \geq \Ll(H_t)$ for all $\nu > 0$.
If we can find conditions on $M, H, K$ such that the capacity-area
inequality holds for these quasi-cylinders,
we could then conclude that both areas must be greater or equal to
$\Ll(H_t)$, and hence, by Lemma~\ref{le:compare}, that
$$
\Ll(K_t) \ge \Ll(H_t).
$$
The above argument shows:

\begin{prop}\label{prop:general} Let $M$ be any symplectic manifold and
$H_{t \in [0,1]}$ a
Hamiltonian generating an isotopy $\phi_t$ from $\1$ to $\phi$. Suppose that
there exists a capacity $c$ such that the following two conditions hold:
\begin{description}
\item[(i)]  $c(H) = \Ll(H_t)$  and

\item[(ii)]   there exists a class $\Cc$ of Hamiltonian isotopies homotopic rel
endpoints to $\phi_{t \in [0,1]}$, which is  such that
the capacity-area inequality holds (with respect to the given capacity $c$)
for all quasi-cylinders
$R_{H,K}(\nu)$ and $R_{K,H}(\nu)$
corresponding to Hamiltonians $K_{t \in [0,1]} \in \Cc$.
\end{description}
Then the length of the path $\phi_t$ is minimal among all paths in $\Cc$.
\end{prop}

Here the class $\Cc$ might be a small neighbourhood of
$\phi_{t\in [0,1]}$, or the set
of all paths from $\1$ to $\phi$
which are homotopic rel endpoints to $\phi_{t\in[0,1]}$.

In \S3  we calculate $c_G(H)$ and $c_{HZ}(H)$ by geometric methods, postponing
to \S4 consideration of the conditions needed on $M$ and
$\Cc$ so that  condition (ii) of Proposition~\ref{prop:general} holds.

\begin{remark}\label{rmk:split}\rm
It is possible that the
Non-Squeezing Theorem holds for all quasi-cylinders.  For example, this would
be the case  if all quasi-cylinders were split, that is
symplectomorphic to a product $M\times D^2$.  When $M$ has
dimension $2$ it follows from the structure theory of ruled surfaces
that any quasi-cylinder splits: see Lemma~\ref{le:prod1}.
Therefore, in this case Theorem~\ref{thm:monotone} is an immediate consequence of
Proposition~\ref{prop:general}.  However, what happens in higher dimensions is still
unclear. Therefore, our arguments are   more indirect and work only under certain
assumptions on the behaviour of holomorphic spheres in $M$.  \end{remark}

\subsection{Properties of quasi-cylinders}

\begin{lemma}\label{le:area}  The area of a quasi-cylinder $Q$ is 
independent of the
choice of the compact set $X$ in Definition~\ref{qc}. 
It is equal to the area
$$
\int_{D_x}\,\Om
$$
 of the disc $D_x = \{x\} \times D$ for any $x\in M$. 
\end{lemma}

\proof{}  First observe that the integral $A_x = \int_{D_x}\Om$ is independent
of the choice of $x$.  This holds because the circles $\{x\}\times \p D$ are
tangent to the characteristic flow of $\Om$ on $\p Q$.  Thus, given any path
$\al$ in $M$, the form $\Om$ restricts to zero on $\al\times \p D$.
We will write $A$ for the common value of the $A_x$.  

When $M$ is compact, we may take $X = M$ and there is nothing further
to prove.  Therefore we now assume that $M$ is noncompact.  

We wish to show that for any $X$ as in Definition~\ref{qc},
$$
\vol\,(X\times D) = A\,\vol\,(X),
$$
that is to say
$$
\int_{X\times D}\,\Om^{n+1} = (n+1)\int_x\om^n\;\int_{D_x}\Om.
$$
But the definition of quasi-cylinder implies that
$$
\Om = (\om\oplus\si)+\tau,
$$
where $\si$ is the standard  area form on $D$, and where $\tau$ is a closed
$2$-form vanishing near the boundary of $X\times D$ and whose restriction to
each fiber $M\times \{pt\}$ vanishes.  Thus $\tau$ has the form
$
\tau = \al\wedge du + \be\wedge dv$, where $\al,\be$ are $1$-forms depending on
the point $p = (x,u,v) \in M\times D$ but with values in $T_x^*M\subset
T_p^*(M\times D)$.  Therefore,
$$
\int_{X\times D}\,\Om^{n+1} = (n+1)\int_{X\times D}\om^n\wedge\si 
+ C\int_{X\times D}\om^{n-1}\wedge \tau^2,
$$
where $C$ is some integer.  Now, the integral of $\tau$ over any element of
$H_2(M;\R)$ $=H_2(M\times D;\R)$ is zero because its restriction to $M\times
\{pt\}$ vanishes.  It is therefore exact, say $\tau = d\rho$.  But then
$\om^{n-1}\wedge \tau^2 = \om^{n-1}\wedge d(\rho\wedge\tau)$, and so, by
Stokes' theorem, the second term in the right hand side of the above equation
must vanish. \QED

The next lemmas establish simple criteria for quasi-cylinders to split.

 \begin{lemma}\label{le:discs}  Let $(M \times U, \tau)$ be a
quasi-cylinder such that the restriction of $\tau$ to each disc 
$pt \times U$ is
non-degenerate.  Then $(M \times U, \tau)$ is split.
\end{lemma} 
\proof{}  By hypothesis $\tau$ has the form $\om + \rho\si +
du\wedge\al + dv\wedge \be$, where the function $\rho $ is $>0$, and $\al =
\al_{u,v}$ and $\be = \be_{u,v}$ are $1$-forms on $M$ for each $u,v\in
U$.  Further, near the boundary $\rho = 1$ and $\al = \be = 0$. Thus
\begin{eqnarray*} 
\tau^{n+1} &=& \rho
(n+1)\om^n\wedge\si -\frac{n(n+1)}2 \om^{n-1} (du\wedge
dv\wedge\al\wedge\be) \\ &=& (n+1)\om^n\wedge\si (\rho + \nu),\;\mbox{
say,}
 \end{eqnarray*} 
where  $\nu > -\rho$. Let 
$$
\tau_t = \om + (1-t+ t\rho )\si+ t(du\wedge\al +
dv\wedge \be).
$$
Then, it is easy to check that
\begin{eqnarray*}
\tau_t^{n+1}& = &(n+1)\om^n\wedge\si (1-t + t\rho + t^2\nu)\\
& > & (n+1)\om^n\wedge\si (1-t + t\rho - t^2\rho) \\
& = & 
(1-t)(1 + t\rho)(n+1)\om^n\wedge\si.
\end{eqnarray*}
Since $(1-t)(1 + t\rho) \ge 0$ when $0 \le t\le 1$, $\tau_t$ is a
symplectic form for all $t$.  Since $\tau_t = \tau_0$ near the boundary,
the usual Moser argument provides an isotopy $h_t$ which is the identity
near the boundary and is such that $h_t^*(\tau_t) = \tau_0$. \QED

\begin{lemma}\label{le:split} When the Hamiltonian $K_t$ is sufficiently close to
$H_t$ in the $C^2$-topology, the quasi-cylinders $R_{H,K}(\nu), R_{K,H}(\nu)$
constructed above   split
for all $\nu$.
\end{lemma}
\proof{}  By symmetry, it suffices to prove this for $R_{H,K}(\nu)$.  
Recall that
the group $\Ham^c(M)$ is locally contractible:  there is a contractible
neighbourhood $\Uu$ of $\1$ which consists of all Hamiltonian diffeomorphisms
$\psi$ whose graph lies close enough to the diagonal $diag$ in $(M\times M,
-\om\oplus \om)$ to correspond to an exact $1$-form in $T^*M$.  (As usual, this
correspondence is induced by a symplectomorphism from a neighbourhood of the
Lagrangian submanifold ${ diag}$ in $(M\times M, -\om \oplus \om)$ to a
neighbourhood of the zero section in the cotangent bundle.)  Thus, if $K_t$ is so
close to $H_t$ that the corresponding paths $\{\psi_t\}, \{\phi_t\}$ satisfy
$\phi_t\circ\psi_t^{-1}\in \Uu$ for all $t$, there is a a canonical choice 
of the
retracting homotopy $f_{u,v}$.  Hence, one can choose the  diffeomorphism 
$\Phi_K:
M\times U_H \to R_{H,K}$ to vary continuously (in the $C^\infty$-topology) with
respect to $K_t$. Observe also, that by choice of $\Phi_K$, 
$$
\Phi_K^*(\Om) = \om \oplus \si
$$
on $R_{H,K}(\nu) - R_{H,K}$.  Therefore, by choosing  $K_t$ sufficiently close to
$H_t$, we may ensure that the form $\Phi_K^*(\Om)$ satisfies the conditions of
Lemma~\ref{le:discs} for all $\nu$. \QED

The results which follow will be needed in \S4.

Let $Q= (M
\times D^2, \Om)$ be a quasi-cylinder of area $A$.  Then there exists a small
neighbourhood of the boundary of the form $M \times C$ for a round annulus
$C$ of radii $r_2 > r_1 >0$,
over which $\Om$ restricts to a split form $\om \oplus (du \wedge dv)$.
Let $C' \subset C$ be the annulus of radii $r_2, r_1 + \frac{r_2 - r_1}{2}$,
and write $Q(\ka), \, \ka>0$ for the quasi-cylinder whose underlying space is
the same as $Q$,
and symplectic form is $\Om_\ka = \Om + \ka\pi^*(\rho)$
where $\pi: Q \to D^2$ is the projection and $\rho$ is a $2$-form supported
in ${\rm Int}\, C'$ which is
everywhere a non-negative multiple of the area form.  We normalise $\rho$
so that $\int \rho = 1$. Then $\area \, Q(\ka) =
A + \ka$. We will denote by $Q^*(\ka)$  the
obvious $S^2$\/-compactification of $Q(\ka)$, which has same area as $Q(\ka)$.
Since $\Om$ is a product in $\pi^{-1}C'$, the forms $\Om_\ka$ are all
symplectic.

 We claim that: 

\begin{lemma}\label{le:qcyl} $Q(\ka)$ is symplectomorphic to a
product for large $\ka$.
\end{lemma}

\proof{}  We may
assume that $Q(\ka)$ is the manifold
$$
Q(\ka) = M \times (U \cup L_\ka) = M \times [0,A+\ka] \times [0,1]
$$
provided with the form $\Om_\ka$ which equals $\Phi^*(\Om)$ in $M \times
U$ and equals the product $\om \oplus \si$ over $L_\ka$, where $\si =
du \wedge dv$. This manifold is embedded as a submanifold of
$V = M \times \R^2$ with the split symplectic form outside $U$.

 Our
aim is to construct a family of disjoint transverse symplectic discs in
$Q(\ka)$ and then apply  Lemma~\ref{le:discs}.  In order that the
diffeomorphism be the identity near the boundary, it is essential
that these discs be flat there, that is they should coincide
with the discs $\{(x,u,v): (u,v) \in  U \cup L_\ka\}$ near the boundary. We will
first construct these discs over $U$ and then extend them over
$L_\ka$.

     For each $x \in M$, let
$D_x$ be the surface obtained by restricting to $M \times U$ all characteristic
lines
of the hypersurfaces $v = $ constant which pass through
the line segment $\{(x,0,v):  0 \le v \le 1\}$.
Because the form
is split near the boundary of $U$, $D_x$ is flat there, and because
$\Om$ restricts to $\om$ on all $M$\/-fibers, these surfaces are everywhere
transverse to the fibers. Hence they are symplectic surfaces transverse to $M$
which foliate $M \times U$, and are standard outside a compact subset of $M$
as well as near the three sides $u=0$, $v=0$ and $v=1$ of $U$. Of course
the map $f_{u,v}:M \to M$ defined by the
identification of $M \times \{(0,v)\}$ with $M \times \{(u,v)\}$ given by
the foliation, is a $\om$\/-symplectic diffeomorphism for all $(u,v)$.
(For $A' > A$ large enough, it is not difficult to see that the isotopy
$f_{A', v \in [0,1]}$ is actually Hamiltonian: the graph of its
generating Hamiltonian is the image by $\psi_T$ of $M \times \{0\}
\times [0,1]$ where $\psi_t$ is the Hamiltonian generated by $F = v$ and
$T$ is large enough so that this image lies entirely outside $U$.)

Let $\be: [A,A+\ka] \to [0, A]$ be a decreasing surjective function which is
constant near the endpoints.  We extend $D_x$ to the disc
$$
D_x' = D_x \cup\{(f_{\be(u),v}(x), u,v): (u,v)\in L_\ka \}.
$$
Clearly, if $\ka$ is sufficiently big, we may choose $\be$ to have so
small derivative that the disc $D_x'$ is symplectic for all $x$.
Further, if $M$ is non-compact, we are only concerned with $x$ in a
compact subset of $M$ since everything is constant near infinity.
Finally, we can use the family of symplectic diffeomorphisms $f_{u,v}$
to define a diffeomorphism $\Psi$ of $M \times [0,A+\ka] \times [0,1]$
which sends the surfaces $D_x$ to flat surfaces. The image of the
symplectic form by $\Psi$ is then still a quasi-cylinder which now
satisfies the condition of Lemma~\ref{le:discs}.
\QED

When $\dim M = 2$, we can
give  a  direct proof that the  manifold $(Q(0),\Om)$ (or equivalently
$(R_{H,K}, \Om)$) is a product.  The reason why this case is
special is that in this dimension we can construct the family  of
disjoint embedded discs or spheres which is needed in
Lemma~\ref{le:discs} as a family of $J$-holomorphic curves. In higher
dimensions there is nothing which forces these curves to be disjoint and
embedded, but in dimension $4$ (the dimension of $M\times S^2$) we can use
positivity of  intersections.

We will need the
following lemma from~\cite{RULED}:

\begin{lemma}\label{le:ruled}  Suppose that $\Om$ is a symplectic form on
$M \times S^2$ which does not vanish on one fiber  $M\times \{z_0\}$ and on
one section $\{x_0\}\times S^2$.  Then, if $\dim M = 2$, $\Om$ is
symplectomorphic to a product form, by a symplectomorphism which is the
identity on $M\times \{z_0\}$ and on $\{x_0\}\times S^2$.
\end{lemma}

\begin{lemma}\label{le:prod1} When $\dim M = 2$ the
manifold $(Q,\Om)$  is a product.
\end{lemma}
\proof{} (Sketch) Since $\Om$ does not vanish on any of the fibers of $Q$,
it suffices by the above lemma to produce one section on which $\Om$ does
not vanish. We do this by showing that for a  generic $\Om$-tame almost
complex structure $J$ on $Q^c$ the class $[pt \times S^2]$ has
an embedded $J$-holomorphic representative.  This is certainly true for
$J$ which are $\Om_{\ka}$-tame, for large $\ka$, and it remains true as $J$
deforms to be $\Om$-tame. For more details, see~\cite{RR,RULED}.\QED

\section{Construction of symplectic embeddings and Hamiltonians}\label{ss:balls}

This section begins by 
constructing functions defined under and over the graph of 
autonomous Hamiltonians $H$ in order to estimate $c_{HZ}(H)$.   
We then
describe situations in which it is possible to calculate
$c_G(H)$  
by embedding balls symplectically under and over the graph of $H$.
(This gives lower bounds of $c(H)$ for any other capacity because, as an
obvious consequence of its definition, $c_G$ is the smallest
capacity.) 
The first case is a local result, where we embed a small ball near an 
extremum of $H$, and
is enough to prove   Theorem~\ref{thm:geodsuff} characterizing geodesics. The
second is a global embedding method which works when $M$ is $2$-dimensional,
and allows us to calculate the Hofer norm of certain elements of $\Ham^c(M)$.

\subsection{Construction of optimal functions}

\begin{prop} \label{prop:HZ-Hamiltonians} Let $M$ be any symplectic manifold
and $H:M \to \R$ be any compactly supported Hamiltonian with no non-constant
closed trajectory in time less than $1$. Then
$$
c_{HZ}(H) \geq \Ll(H).
$$
\end{prop}

\proof{} Fix some small $\nu > 0$. We must show that $c_{HZ}(R_H^-(\nu/2)) \geq 
\max H$, where $H$ is rescaled so that $\min H = 0$. Set $\max H = m$ and consider
the following space
$$
S_{H, \nu/4} = \{ (x,\rho,t) \in M \times D^2(m+\nu/4) \; \mid \; 
0 \leq \rho \leq H(x) + \nu/4 \}.
$$  
Here $(\rho,t)$ are the action-angle coordinates of the disc, and are
related to polar coordinates $(r,\theta)$ by
$$
\rho = \pi r^2/2,\quad \theta = 2\pi t.
$$
It is very easy to see that $S_{H, \nu/4}$ can be embedded in $R_H^-(\nu/2)$ by
a symplectic map which moves points only in the $\R^2$-factor, that is to say by
the restriction to $S_{H, \nu/4}$ of a map of the form $\id \times \phi: 
M \times \R^2 \to M \times \R^2$. Thus it suffices to show that $c_{HZ}(S_{H, \nu/4})
\geq m$. Define the function $K: S_{H, \nu/4} \to \R$ by
$$
K(x,\rho,t) = -H(x)  + m + \rho = -H(x)  + m + \pi r^2.
$$
 Then clearly $K$ is a smooth positive function which has no non-constant
closed trajectory in time less than $1$, and is constant on the boundary of 
$S_{H, \nu/4}$. But $\| K \| = m + \nu/4$ which shows that $c_{HZ}(S_{H, \nu/4}) 
\geq m$, at least when $M$ is compact. If $M$ is non-compact, and since 
${\rm supp} H$ is compact, the above function $K$ is not equal to $m + \nu/4$ outside
a compact subset of  $S_{H, \nu/4}$. This is because $K$ is non-constant on
$(M - {\rm supp} H) \times D^2(\nu/4)$. But we can smooth out $K$ to the constant
value $m + \nu/4$ on that subset. Clearly, this can be done without creating a
closed trajectory in time less than $1$ as soon as $\nu$ is small enough. 
\QED

\subsection{Local embeddings}

Recall that a path $H_{t\in [0,1]}$ is called quasi-autonomous if 
it has a fixed minimum and a fixed maximum.

\begin{lemma}\label{le:C2} There is a $C^2$-neighbourhood $\Uu$ of 
$0$ in the space of
Hamiltonians $M \times [0,1] \to \R$ such that 
$$
c_G(H) = \Ll(H_t)
$$
for all quasi-autonomous paths $H_{t\in [0,1]} $ in $\Uu$.
\end{lemma}

\proof{}  Let $P$ and $p$ be  fixed extrema of $H_t$. After reparametrization,
we may assume that $\min H_t = 0$.  As before, let $\mu_H(t) = \max_x H_t(x)$
and set 
$$
m = \int_0^1\mu_H(t)\,dt =\Ll(H_t).
$$
  We will show how to embed the ball
$B^{2n+2}(m - \eps)$ in $R_{H}^-$ for all $\eps > 0$. The argument
for $R_{H}^+$ is similar.

We will construct a map $\tilde f: B^{2n+2}(m-\eps) \to R_{H}^-$ such
that  the
following diagram commutes:
$$
\begin{array}{cccc}B^{2n+2}(m - \eps) &\quad \stackrel{\tilde f}{\to} \quad
& R_{H}^\pm &\subset\; M\times \R^2\\ \\
         \pi \downarrow & & \downarrow \pi'  & \\ \\
B^{2n}(m - \eps) & \stackrel{f}{\to} & M &
\end{array}
$$
where the vertical arrows are the obvious projections.  
(Such fibered embeddings
were considered also in~\cite{LALMCD}.) 
Observe that if $x \in \p B^{2n}((1-c)(m-\eps))$ then the fiber
$\pi^{-1}(x)$ is the $2$-disc $B^2(c(m-\eps))$, which sits inside 
$B^2(m-\eps)$.  Hence we may construct $\tilde f$ as the product
of a symplectic embedding $f$ with a suitable  area-preserving
embedding $\psi$ of the $2$-disc $B^2(m-\eps)$ into $\R^2$.

To construct $f$,  observe first that if $H_t$ is $C^2$-small,
there clearly is a symplectic embedding $f:B^{2n}(m-\eps) \to  M$ 
which
takes $0$ to the fixed maximum $P$ and is such that, for each $c\in [0,1]$,
$$
f(B^{2n}((1-c)(m-\eps)))\subset H_t^{-1}[c\mu_H(t), \mu_H(t)],\;\;\mbox{for
all}\;\;t\in[0,1]. 
$$
Then,  if $x \in \p B^{2n}((1-c)(m-\eps))$, the fiber
$(\pi')^{-1}(f(x))$ contains the set
$$
U_c = \{(s,t)\in \R\times [0,1]: 0 \le s \le c\mu_H(t)\}.
$$
Since  the area of $U_c$ is greater or equal to $cm$ and $\eps > 0$, 
it is easy to see that
there is an area-preserving embedding $\psi$ of the $2$-disc
$B^2(m - \eps)$ into $U_1$ which takes $B^2(c(m - \eps)) $ into $U_c$
for each $c\in [0,1]$.
This completes the proof.\QED

\begin{remark}\rm Note that because the Gromov capacity $c_G$ is the smallest
capacity, the above lemma proves that $c_{HZ}(H) \geq \Ll(H_t)$ too.
\end{remark}

   Here is an immediate corollary.

\begin{cor}\label{co:lengmin2} 
Let $\{H_t\}$ be a Hamiltonian generating a regular path,  which has
at least one fixed minimum and one fixed maximum at each moment. Then,
there exists $\eps_0 >0$ such that for all $ 0 < \eps \le\eps_0$,
$$
c_G(\eps H_{a + t\eps}) = \Ll(\eps H_{a+t\eps}).
$$
\end{cor}

Here, of course, $H_{a+t\eps}$ is short for $H_{a+t\eps}, t\in [0,1]$.
We can now prove:

\begin{cor}  Theorem~\ref{thm:geodsuff} holds.
\end{cor}

\proof{}  Suppose that $H_{t\in [0,1]}$ is a Hamiltonian which has a fixed maximum
and a fixed minimum at each moment, and let   $\eps_0$ be as above.  Then we
claim that, whenever  $0 < b-a \le\eps_0$, the path $\phi_{t\in[a,b]}$ is a stable
geodesic. Since its generating Hamiltonian is just $(b-a)H_{a+t(b-a)}$, the above
corollary shows that  condition (i) of Proposition~\ref{prop:general} holds for
$c = c_G$.   If we take
$\Cc$ to be a sufficiently small neighbourhood of $\phi_{t\in[a,b]}$, then condition
(ii) holds also by Lemma~\ref{le:split} and Theorem~\ref{thm:ns}.  Therefore, the
result follows from Proposition 2.2.\QED

\subsection{Global embeddings when $\dim\, M = 2$}

We  discuss the construction of symplectic embeddings for Hamiltonians defined
on surfaces. Although the method of estimating $c_{HZ}(H)$ given in \S3.1
is easier, one sometimes obtains sharper estimates by using embedded balls:
see Remark~\ref{rmk:cap}.
However, this section is not used in the applications in \S5,
except for the $S^2$ case of Theorem~\ref{thm:last}.

\subsection*{Trapezoids}

It is often convenient to think of embedding linear shapes such as
trapezoids instead of  balls.  The following lemma shows
that the embedding problems for these different shapes are
essentially the same.

As above, we write $B^{2n}(c)$ for the ball with capacity $c$.
Let $T^{2n+2}(a)$ be the trapezoid:
$$
T^{2n+2}(a) = \{(x,u,v)\in \R^{2n+2}: 0 \leq u\leq a, 0\leq v\leq 1
, x\in B^{2n}(a-u))\}.
$$
We think of $T^{2n+2}(a)$ as fibered over the $(u,v)$ rectangle $R(a)$,
with  suitable balls as fibers.  Similarly, $B^{2n+2}(c)$ fibers over
$B^2(c)$ with balls of varying size as fibers.  The embeddings which we now
construct  will preserve this fibered structure.

\begin{lemma}\label{le:trap}  For all $\eps > 0$,
\begin{description}
\item[(i)]
  $B^{2n+2}(a)$ embeds symplectically in $T^{2n+2}(a+ \eps)$,

\item[(ii)]  $T^{2n+2}(a)$ embeds symplectically in $B^{2n+2}(a + \eps)$.
\end{description}
\end{lemma}

\proof{}  $B^{2n+2}(a)$ fibers over the disc $B^2(a)$.  Let $h_B:
B^2(a) \to \R$ be the function which at the point $(u,v)$ equals the
capacity of the fiber at that point.  So $h_B(u,v) = a - \pi(u^2 +
v^2)$.    Let $h_T$ be the similar function for $T^{2n+2}(a)$ (defined
over the rectangle  $R(a)$).  To define an embedding  $B^{2n+2}(a) \to
T^{2n+2}(a+ \eps)$ which is the  identity on the fibers, it clearly
suffices to find an area preserving map $f :
B^2(a) \to R(a+\eps)$ such that
$$
h_T(f(u,v)) \geq h_B(u,v).
$$
It is not hard to see that such a map $f$ exists.   For example, one
can construct $f$ as follows.  It should
take $ p = (-\sqrt{a/\pi},0)$ to the vertex $(0,0)$
of $R(a + \eps)$.  Further, if
$\ga_\mu$ is a family  of disjoint loops based at $p$ which go
close to the ray $y = 0$ then round the circle centered at the origin
with radius $\mu$  and
then back again close to the ray, $f$ should take this
family to an appropriate family of loops in the rectangle.  These
loops exist because,
 for all $\la$,
$$
area\,\{(u,v): h_B(u,v) \geq \la\}\; \le\; area\,\{ (u,v): h_R(u,v) \ge
\la\}.
$$
One also uses the fact  that these sets are connected.  This shows (i).

The proof of (ii) is similar. \QED

We now apply these ideas to calculating $c_G(H)$.  We will give a prototype result
here, with only a partial proof.  Full details of the proof and extensions of this
result may be found in~\cite{RICH}.

\begin{prop}\label{prop:ball4}  Let $(M,\om)$ have dimension $2$,  and
let $H$ be an autonomous Hamiltonian with time-$1$ map $\phi$.
Suppose that $H$ has no non-constant
$t$-periodic orbits for $0 < t < 1$. Then
$$
c_G(H) = \|H\|.
$$
\end{prop}

\proof{}  We will prove this only in the special case when there is a path $\be$ in
$M$ from an absolute maximum $P$  to an absolute minimum $p$ of $H$ along
which $H$ decreases. In other words, we assume that the gradient flow of $H$ with
respect to a generic metric on $M$ has a trajectory going from $P$ to $p$.
  We will renormalise
$H$ so that $H(P) = m, H(p) = 0$.  Thus $\|H\| = m$.

We will show that, given any $\eps > 0$, we can
 embed the
trapezoid  $T = T^4(m - \eps)$ under the graph of $H$ in $X = M \times
[0,m] \times [0, 1]$.  
By Lemma~\ref{le:trap},  this will imply that $c_G(R_H^-) = m$.  By replacing $H$
with $m - H$, we obtain $c_G(R_H^+) = m$.  Hence
$$
c_G(H) = \min\{c_G(R_H^-), c_G(R_H^+)\} = m
$$
as required.

 Let  $\xi_H$ be the Hamiltonian vector field on $M$ corresponding to $H$.
By slightly perturbing $H$, we may assume that $\xi_H$ has a finite number of
critical points. Let
 $M'$ be $ M$ minus these critical points. Then the given area form
$\si$ on  $M$ may be written as $\al \wedge dH$ on $M'$ where
$\al$ is a $1$-form such that $\al(\xi_H) = 1$.  Observe that, because
$H$ has no periodic orbits of period $< 1$, each level set of $H$ has
$\xi_H$-length at least $1$.  Thus, the integral of $\al$ round each
level surface of $H$ is at least $1$ and so the area of $\Sigma$ is at
least $m$.

We will construct a fibered embedding of ${\rm Int}\, T$ into $X$ which
covers an embedding $f: {\rm Int}\,R(a) \to M'$.  We define $f$ as
follows.   Parametrize the path $\be$ from $P$ to $p$ so that
$$
H(\be(t)) = m - t,\quad t\in [0,m],
$$
and define $f$ along the axis $v = 0$ by
$$
f(u,0) = \be(\eps/2 + u),\quad u\in [0,m-\eps].
$$
Then define $f$ so that it takes each   arc $u = c$ 
into the level set $H = m -\eps/2- c$ in such a way that $\frac
{\p} {\p v}$   goes to $\xi_H$.  This is an embedding because the arcs $H
= const$ all have $\xi_H$-length at least $1$.  Moreover, it preserves
the $u,v$ area, since  $dH$ pulls back to $-du$.

Observe that the fiber in $X$ over the point $f(u,v)$ is the rectangle
$R_c = [0,c]\times [0,1]$ where $c = {m-u-\eps/2}$.  It is
not hard to see that there is an area preserving map
$$
g:B^2(m-\eps) \to R_{m-\eps/2}
$$
which takes $B^2(c-\eps/2)$ into $R_{c } $ for all $c$.  Then the product
$$
g\times f: T\to X
$$
is the desired embedding.
\QED

\section{Generalised capacity-area inequalities}\label{ss:global-ns}

We now investigate conditions on $M$ and classes $\Cc$ so that capacity-area
inequalities hold for spaces of the form $R_{H,K}(\nu)$ where $H,K \in \Cc$.
When combined  with Proposition~\ref{prop:general} this will prove
Theorem~\ref{thm:monotone}.  

We begin by considering the Gromov capacity.
The proof in \cite{LALMCD} of the Non-Squeezing Theorem in all manifolds
does not extend to quasi-cylinders $Q$, the main difficulty being
that, if one considers $Q$ lying inside $M \times \R^2$, it is not clear that
one can disjoin $Q$ from itself with energy equal to $\area \,(Q)$. But,
fortunately,  a holomorphic argument, similar to the one contained in
Proposition~3.1 of
\cite{LALMCD}, can be applied here (see
Proposition~\ref{prop:deformation1} below).  We begin with
a different argument, where we first extend the quasi-cylinder and then
trivialise the extension.

\begin{prop}\label{prop:quasi-ns} Let $Q = (M \times D^2, \Om)$ be a
quasi-cylinder. Then the Non-Squeezing Theorem holds whenever $M$ has
bounded geometry at infinity and satisfies one of the  following conditions
\begin{description}
\item[(i)] $\dim(M) \le 4$;
\item[(ii)] $(M \times S^2(a), \om \oplus \si)$ is weakly monotone for all
$a>0$.
\item[(iii)]  there are no spherical
homology classes $A\in H_2(M)$ such that
$$
\om(A) > 0,\quad 2-n \le c_1(A) \le 0,
$$
\end{description}
\end{prop}

\begin{remark}\label{rmk:condit}\rm  The condition in (iii) above is clearly
satisfied when $M$ is semi-monotone.   However, it is also satisfied when the
minimal Chern number (the smallest positive value of $c_1$ on spherical
homology classes) is $\ge n-2$, for instance when $M = S^2\times S^2\times S^2$,
and for various other manifolds. \end{remark}

\proof{}
We will begin with the proof of the sufficiency of condition (ii).  Since
manifolds of dimension $\le 6$ are automatically weakly monotone,  this will
imply the sufficiency of condition (i).
\MS

\NI
{\bf Proof of the sufficiency of condition (ii) of
Proposition~\ref{prop:quasi-ns}}. \quad
\smallskip
 
 Let us assume, by contradiction, that $Q$ contains a ball
of capacity $c > A$, and consider its split expansion $Q(\ka)$ 
defined in Lemma~\ref{le:qcyl}. Then,
by construction,
$Q(\ka)$ contains a disjointly embedded ball $B^{2n+2}(c)$ and
cylinder $M\times B^2(\ka)$, where $c + \ka > A + \ka$.
We now use the theory of $J$-holomorphic curves to show that this is
impossible.  We assume that the reader is familiar with the basics of
this theory, as explained in \cite[\S3]{LALMCD} or~\cite[Chapter~1]{DUSAL}, for
example.  This theory tells us that, given any $\Om$-tame almost
complex structure $J$ on the compactification $Q^*(\ka) = M\times S^2(A+\ka)$ and
 any point $x \in Q^*(\ka)$, there is a
$J$-holomorphic curve (or cusp-curve) $C$ through $x$, which
represents the homology class of $pt \times S^2$ and so satisfies
$$
\int_C\Om = A + \ka.
$$
(It is at this point that we use the hypothesis that $M\times S^2(A+\ka)$ is
weakly monotone.  For general manifolds, it is not known whether the above
statement is true for any $J$.)

 Let us apply this when $J$ restricts to the
standard complex structure $J_0$ on the ball $B$ and  to a product
structure $J_{spl}$ on the cylinder ${\rm Cyl}$, and let us take the point
$x$ to be the center of the ball.  Then, as usual, it follows by Gromov's
monotonicity argument that $$
\int_{C\cap B} \Om \ge c.
$$
We claim that
$$
\int_{C\cap {\rm Cyl}} \Om \ge \ka.
$$
 To  see this, observe that the projections
$\pi_D: {\rm Cyl} \to D = B^2(\ka)$ and $\pi_M:{\rm Cyl} \to M$ are
both holomorphic.  Thus, if $\si$ denotes the area form on $D$,
\begin{eqnarray*}
 \int_{C\cap {\rm Cyl}} \Om & = & \int_{\pi_D(C\cap{\rm Cyl})} \si
                                         + \int_{\pi_M(C\cap{\rm Cyl})}\om \\
& \ge & \int_{\pi_D(C\cap{\rm Cyl})} \si\\
& \ge &  \ka,
\end{eqnarray*}
where the first inequality holds because $\pi_M(C\cap{\rm Cyl})$ is a
holomorphic
curve in $M$ and the second holds because  $\pi_D: C\cap{\rm Cyl} \to
D$ is surjective. ($C$ must intersect every
fiber $M \times pt$ of $X \to S^2$  for topological reasons.)

Since $J$ is $\Om$-tame, the integral of $\Om$ over the piece of $C$
which is outside $B\cup{\rm Cyl}$ is
positive.  Therefore,
$$
\int_C\Om >\int_{C\cap B} \Om +  \int_{C\cap {\rm Cyl}} \Om \ge c +
\ka > A+\ka = \int_C\Om,
$$
which is impossible.
\MS

\NI
{\bf Proof of the sufficiency of condition (iii) of
Proposition~\ref{prop:quasi-ns}}. \quad
Here we need  to control the degeneracies of the pseudo-holomorphic
sphere. We do this by considering a special (not necessarily generic)
path of almost
complex structures. Let $J_t$ be
a path of $\Om$-tame almost complex
structures on  $Q^*(\ka)$ such that $J_0$ is the split stucture and
$J_1$ is standard with
respect to the ball and to the cylinder, and such that one
fiber $M\times pt$ is $J_t$-holomorphic for all $t \in [0,1]$.  
(It is easy to find
such a path.)  It suffices to show that there is a
$J_1$-holomorphic curve in  class $\{pt\} \times S^2$ through the center $p$
of $B$.
We first show that along such a path of almost complex structures,
only very special degeneracies can occur.
 Indeed, let $C$ be
any  $J_t$-holomorphic cusp-curve
$$
C = C_1\cup\dots\cup C_k
$$
in  class $[\{pt\} \times S^2]$.
We may decompose the homology class of $C_i$ as
$$
[C_i] = p_i[pt\times S^2] + Z_i,\quad \mbox{where}\quad Z_i\in H_2(M).
$$
Since every intersection of $C_i$ with the $J_t$-holomorphic fiber
$M\times pt$ contributes positively to the
intersection number $p_i= [C_i]\cdot [M]$, we must have $p_i \ge 0$ for
all $i$.  Hence $p_1 = 1$ and $p_i = 0, i \ge 2$, so
$[C_i] \in H_2(M)$ for $i \ge 2$.  Further, because $J_t$ is
$\Om$-tame,
$$
0 < \int_{C_i}\Om = \int_{C_i}\om.
$$
for all $i \geq 2$.

We claim that a similar statement holds for any path $J_t'$  which is
sufficiently close to $J_t$.  For otherwise, there would be a sequence of
cusp-curves $C^\nu = C_1^\nu\cup\dots\cup C_{k_\nu}^\nu$ converging to a
$J_t$-holomorphic cusp-curve $C = C_1\cup\dots\cup C_k$ with $[C-i]\in
H_2(M)$ for $i\ge 2$.   But the compactness theorem implies that (after taking
an appropriate subsequence) the components of $C^\nu$ converge to unions of
components of $C$.   Hence for large $\nu$, all but one of the classes
$[C_i^\nu]$ must belong to $H_2(M)$.  Therefore, we may assume that the path
$J_t'$ is a generic path which joins the regular element $J_0$ to a point $J_1'$
which is very close to $J_1$.  Then the hypothesis (iii) on $M$ implies that the
only classes $[C_i]\in H_2(M)$ which have $J_t'$-holomorphic representatives
in $Q^*(\ka)$ for some $t$ are those with $c_1\ge 0$ since these are the ones for
which the corresponding moduli space has non-negative dimension.  (The
number $4$ rather than $3$ was used in (iii) to compensate for the fact that
$\dim\, Q^*(\ka) = \dim\,M + 2$.)

Thus none of the $J_t'$-holomorphic cusp-curves contain multiply-covered
components of negative Chern number.  Therefore, by the remarks
in~\cite[\S3]{LALMCD}, there is a $J'_1$-curve in class $[pt \times S^2]$
passing through any generic point arbitrarily close to $p$. Using once again the
compactness theorem, we conclude that there must be a $J_1$-curve in class
$[pt \times S^2]$ passing through $p$.
\QED
\MS

  Our last global Non-Squeezing Theorem holds
for rational manifolds. Let $r$ be the index of rationality of $M$.
For the convenience of the reader we restate
Proposition~\ref{prop:deformation0}.

\begin{prop}\label{prop:deformation1}
Let $\Om_{s \in [0,1]})$ be a smooth
$1$\/-parameter family of symplectic forms on $M\times D$ such that
each manifold  $Q_s = (M \times D, \Om_s)$ is a quasi-cylinder. Assume
that $Q_0$ is split.
\begin{description}
\item[(i)]   If $\area \, (Q_s) < 2r$ for all $s$,
then the non-squeezing theorem holds for $Q_1$;

\item[(ii)]
 If there is a smooth symplectic isotopy
$$
g_s: B^{2n+2}(c) \to Q_s, \quad 0 \le s \le 1
$$
which at time $s=0$ has the form
$$
B^{2n+2}(c) \hookrightarrow B^{2n}(c)\times D(c)\stackrel{f\times g}
{\longrightarrow} M\times D
$$
where $c \ge 2r$, and if
$$
\area \,Q_s < c + 2r, \quad 0 \le s \le 1,
$$
then $\area \, Q_1 \ge c$.
\end{description}
\end{prop}
\proof{}  We begin with the proof of (ii).
 The idea
is to use the embedded ball in the quasi-cylinders in order to
control the area of the components of any possible cusp-curve. Actually,
we will combine the control given by these embedded balls to the one
provided by the positivity of intersection with a given fixed $J$-invariant
fiber.

Suppose, by contradiction, that $\area \, Q_1 < c$.
Denote by $Q^*_s$ the $S^2$-compact\-ification of
$Q_s$ obtained by attaching a cylinder $(M \times D, \om \oplus \si)$
of area $a$ to the boundary of $Q_s$, where $a$ is small enough so that
$A_s = \area \,(Q^*_s) < c + 2r$ and $A_1 < c$. We still denote by
$\Om_s$ the form
on $Q^*_s$. Let $U$ be a small neighbourhood of the point $p_{\infty}$, the
center of $D$, sitting in the interior of $D$.
It is enough to produce a $J$-$A$\/-curve
passing through $p=g_1(0)$, where $J$ extends the standard structure on
$g_1(J_0)$ and is split on the open subset $W= M \times U$.
 Let $J_s, \, 0 \le s \le 1$,
be a smooth (not necessarily generic) one-parameter family of
almost complex structures tamed by $\Om_s$, beginning with the split (regular)
structure and ending at $J_1 = J$, such that $J_s$
be split on $W$ and extends for each $s$ the structure $g_s(J_0)$.
Then any $J_s$-$A$\/-cusp-curve passing through $p_s = g_s(0)$ must be of type
$$
A = (A-B) + B  \quad \; \mbox{or} \quad \; A = (A-B) + B_1 + B_2
$$
for some classes $B, B_1, B_2 \in H^2(M)$.
Indeed, the positivity of intersection
implies as before that $A = (A-B) + \sum_i B_i$, but the component that passes
by $p_s$ must have area at least $c$. There remains then only a total area
less than $2r$ for all other components, which easily implies the above
cusp-curve decomposition. Because $A-B$ is a primitive class, it cannot be
a multiple covering. Consider a $B$\/-component: if it passes through $p_s$,
and is a multiple covering, its area would be at least $2c$, and therefore
would be greater than $\area \,Q^*_s$, since $c \ge 2r$ and $\area \,Q^*_s
< c + 2r$. If it does not pass through $p_s $, its area is smaller than $2r$ and
therefore cannot be a multiple covering either.
Hence no component can be a multiple covering. As above, we conclude
that there are $J$-$A$\/-curves passing through $p$.

The proof of (i) is easier and left to the reader.

\bigskip

   We now turn to $c_{HZ}$-area inequalities for quasi-cylinders.

\begin{prop} \label{prop:HZ-qc} Let $M$ be weakly exact. 
Then the $c_{HZ}$-area inequality
holds for all quasi-cylinders $Q = (M \times D^2, \Om)$.
\end{prop}

\proof{} Assume by contradiction that $c_{HZ}(Q) > \area (Q)$, and let 
$H:Q \to [0, c]$ be such that $H = c > \area (Q)$ on $\p Q$ and outside a
compact subset, vanishes somewhere inside $Q$ and has no non-constant
closed trajectory in time $< 1$. Let $\eps$ be smaller than 
$c - \area (Q)$. Then $Q$ sits in the interior of the cylinder
$Q(\ka + \eps)$ (here $Q(\ka)$ is defined as in Lemma~\ref{le:qcyl}
and the $\eps$-part 
ensures that $Q$ is in the interior of $Q(\ka + \eps)$). Then extend
$H$ to a Hamiltonian $\tilde{H}: Q(\ka + \eps) \to [0, \infty)$
which equals $c + \ka$ on $\p Q(\ka + \eps)$. Because $Q(\ka + \eps) - Q$
is split and has area $\ka + \eps$, one can choose such an extension
so that $\tilde H$ has no non-constant closed trajectory in time $< 1$.
But then $c_{HZ}(Q(\ka + \eps)) \geq c + \ka > \area (Q(\ka + \eps))$, which
contradicts Theorem~\ref{thm:HZ-ineq}.
\QED

\begin{remark}\label{rmk:capsph}\rm
Observe that because $\ka$ may  be arbitrarily large here, 
we cannot say anything
about manifolds with finite index of rationality unless $M = S^2$.  
In this case, by
Lemma~\ref{le:prod1}, the above result holds for 
for all quasi-cylinders of area less than or equal to the
area of $S^2$ because these cylinders all split.
\end{remark}

\section{Applications to Hofer's geometry}\label{ss:length-min}

We now apply the results of \S\ref{ss:balls} and \S\ref{ss:global-ns}
to derive the theorems on length-minimizing paths in Hofer's metric.
Our tool to do this is Proposition~\ref{prop:general}. We first discuss general
results with fixed endpoints (global with respect to time),
and then give results which hold at each moment (local with respect
to time). We then turn to specific results, linear rigidity, and
a brief discussion of Hofer's diameter.
Finally, we prove  the local flatness result  Proposition~\ref{prop:BP}.

\subsection{Length-minimizing paths: general results}

To begin, here is a slightly more general version of Theorem~\ref{thm:monotone}.

\begin{theorem}\label{thm:monotone1} 
Let $M$ have dimension $\le 4$ or be semi-monotone (or, more generally,
satisfy condition (iii) in Proposition~\ref{prop:quasi-ns}).  Assume further in 
the
non-compact case that it has bounded geometry at infinity. Let $H_{t \in
[0,1]}$ be any path with $$
c_G(H) = \Ll(H_t).
$$
Then $H_{t \in [0,1]}$ is length-minimizing amongst all paths homotopic
(with fixed endpoints) to $H_{t \in [0,1]}$. The same conclusion holds
if $c_{HZ}(H) = \Ll(H_t)$ when $M$ is weakly exact.
\end{theorem}

  This is an obvious consequence of Propositions~\ref{prop:general},
\ref{prop:quasi-ns} and \ref{prop:HZ-qc}.
  The next result talks about standard balls in
$R_H^\pm(\nu/2)$.  These are balls which are embedded by a map
which is the restriction to $B^{2n+2}(r)$ of a product embedding:
$$
 B^{2n}(r)\times B^2(r)\to M\times U_H\subset M\times \R^2.
$$

\MS
\begin{theorem}\label{thm:rational}
Let $M$ be a rational symplectic manifold of rationality index
$r$
with bounded geometry at infinity. Let $H_{t \in [0,1]}$ be
any path with
$$
c_G(H) = \Ll(H_t).
$$
Then:
\begin{description}
\item[(i)] In the case $\Ll(H_t) < 2r$, $H_{t \in [0,1]}$ is length-minimizing
amongst all paths which are homotopic (rel endpoints) to $H_{t \in [0,1]}$
through paths of lengths $< 2r$.

\item[(ii)] In the case $\Ll(H_t)\ge 2r$,
assume that, for all small $\eps > 0$, there exist an embedded ball $B$
of capacity $\ge \Ll(H_t) - \eps$ in $R_H^+$ and a number $0 < \nu < 2r$
such that $B$ is
isotopic inside $R_H^+(\nu/2)$ to a standard ball in $R_H^+(\nu/2)$.
Suppose that the same holds for $R_H^-$. Then
$H_{t \in [0,1]}$ is length-minimizing
amongst all paths which are homotopic (rel endpoints) to $H_{t \in [0,1]}$
through paths of lengths $< \Ll(H_t) + 2r$.
\end{description}
\end{theorem}

\proof{} (i) $\;$ Let $\Ll(H) =\Ll(H_t) < 2r$ and $\psi_{t,s \in [0,1]}$
be a homotopy of paths generated by $K^s, s\in[0,1]$ of length
$$
\Ll(K^s) <  2r.
$$
By Proposition~\ref{prop:general}, it is enough to show that the Non-Squeezing
Theorem holds for the quasi-cylinders $R_{H,K^{1}}(\nu),
R_{K^{1},H}(\nu)$, at least when $\nu$ is sufficiently small.
We prove this using Proposition~\ref{prop:deformation1} as follows.

Suppose, by contradiction, that for all small $\nu$, there is an embedded ball
$B_{\nu}$ of capacity
$c_{\nu}$ in say $R_{H,K^{1}}(\nu)$ with
$$
\area \,(R_{H,K^{1}}(\nu)) < c_{\nu}.
$$
Choose $\nu$ small enough so that $\|K^s\| + 2\nu <  2r$ for all $s$, and
fix once and for all that value $\nu$.
Set $c = c_{\nu}$, $Q_s = (R_{H,K^{s}}(\nu), \Om)$ and
$Q'_s = (R_{K^{s},H}(\nu), \Om)$. Up to a smoothly varying symplectomorphism,
we have:
$$
Q_s = (M \times D^2, \Om_s) \quad Q'_s = (M \times D^2, \Om'_s).
$$
Since $\area \, Q_1 < c$, Proposition~\ref{prop:deformation1} (i)
implies that there exists some
$s$ such that $\area \,Q_{s} \ge  2r$. Let $\bar s$ be the smallest value
of $s$ for which such an inequality holds for either $Q_s$ or $Q'_s$.
Suppose that this happens first for $Q$ say. This means that
$$
\frac {1}{\vol\, M}(\vol(R_H^-) + \vol(R_{K^{\bar s}}^+)) + \nu \ge  2r
$$
or equivalently:
\begin{eqnarray}
\la \Ll(H) + \mu \Ll(K^{\bar s}) + \nu \ge 2r     \label{ineq:lambda-mu}
\end{eqnarray}
for some $\la, \mu \in (0,1)$. Because $\bar s$ is the smallest such value,
$$
\area \,(Q'_s) < 2r
$$
for all $s < \bar s$. But, since $R_H^+$ embeds in $Q'_s$, there exist
embedded balls of capacity arbitrarily close to $\Ll(H)$ inside $Q'_{\bar s}$.
Therefore, using again Proposition~\ref{prop:deformation1} (i), we must have
$$
\area \,(Q'_{\bar s}) \ge \Ll(H),
$$
or equivalently:
\begin{eqnarray}
(1-\la) \Ll(H) + (1-\mu) \Ll(K^{\bar s}) + \nu \ge \Ll(H).
\label{ineq:mu-lambda}
\end{eqnarray}
Then the sum of inequalities~(\ref{ineq:lambda-mu}) and~(\ref{ineq:mu-lambda})
gives
\begin{eqnarray*}
\Ll(H) + \Ll(K^{\bar s}) + 2 \nu & \ge & \Ll(H) + 2r \\
        \Ll(K^{\bar s}) + 2 \nu & \ge &  2r
\end{eqnarray*}
a contradiction.

\medskip

(ii)
 This is a consequence of Proposition~\ref{prop:deformation0} (ii).
The proof, similar to the proof of (i), is left to the reader.

\QED

\MS
   The results of \S\ref{ss:balls} and \S\ref{ss:global-ns} also lead to the
following theorem stating that geodesics minimize the length at each moment
in a strong sense.

\begin{theorem}
Let $M$ be a symplectic manifold with
bounded geometry at infinity, and let $\phi_{t \in [0,1]}$ be
a geodesic.
\begin{description}
\item[(i)]
Assume that $M$ is semi-monotone. Then  $\phi_{t \in [0,1]}$
minimizes length at each moment amongst all homotopic paths. More
precisely,   each $s\in [0,1]$ has a closed connected neighbourhood $\Nn(s)$
such that   the path $\phi_{t
\in \Nn(s)}$ is length-minimizing amongst all homotopic paths.
\item[(ii)] Assume that M is rational. Then
 each $s\in [0,1]$  has a closed connected neighbourhood $\Nn(s)$  such that
for all closed subintervals $\Nn$ containing $s$  
the path $\phi_{t \in \Nn}$ is  length-minimizing amongst all
paths which are homotopic to $\phi_{t \in \Nn}$ through paths of lengths
$< 2r$.
\end{description}
\end{theorem}

  The statement (i) is a consequence of Theorem~\ref{thm:monotone} and
Corollary~\ref{co:lengmin2}.
The second one is a consequence of Theorem~\ref{thm:rational} (i) and
Corollary~\ref{co:lengmin2}.

\subsection{Length-minimizing paths: specific results}

We now combine the results of the previous sections to estimate the
Hofer norm $\|\phi\|$ of certain elements $\phi\in \Ham^c(M)$. We begin with
the proof of: 

\begin{theorem}~\label{thm:last}
Let $M$ be a weakly exact  manifold
 or any surface, and let $H$ be
an autonomous Hamiltonian with time-$1$ map $\phi$. Suppose that
$H$ has no non-constant closed trajectory in time less than $1$.
Then  
 $\phi_{t\in[0,1]}$ is length-minimizing in its homotopy
class. 
\end{theorem}

\proof{} This is an immediate consequence of Theorem~\ref{thm:monotone}
and Proposition~\ref{prop:ball4} for surfaces and of 
Theorem~\ref{thm:monotone}
and Proposition~\ref{prop:HZ-Hamiltonians} for the weakly exact case.
\QED

\bigskip
\noindent
{\bf Proof of Corollary~\ref{cor:HS-absolute}}. $\;$
When $M$ is a surface and has genus $> 0$, the group $\Ham(M)$ is contractible,
and when $M$ is exact and convex, $\Ham^c(M)$ has
$r_1 = \infty$. Thus in these cases the result follows 
from Theorem~\ref{thm:last}.

      However, $\Ham(S^2)$ deformation retracts to ${\rm SO}\,(3)$.  In
particular, its fundamental group has order $2$, and is generated by
a full rotation around a fixed axis.  Thus the isotopies from
$\1$ to $\phi$ divide into two homotopy classes: those which are
 homotopic to
the flow $\phi_{t \in [0,1]}$ of $H$, and those which are not. 
Theorem~\ref{thm:monotone} and
Proposition~\ref{prop:ball4}
tell us that all isotopies in the former class have length $\ge \|H\|$.
(See also Remark~\ref{rmk:capsph}.)
It therefore remains to consider the length of isotopies in the
latter class.  But, it follows from Lemma~\ref{le:turn1}
below that, for such isotopies
$\psi_t$,
$$
\Ll(\phi_{t \in [0,1]}) + \Ll(\psi_{t \in [0,1]}) \ge \area\, S^2.
$$
Therefore, if $\Ll(\phi_{t \in [0,1]}) = \|H\|$ is no more than half the area,
it has to be minimal.
\QED

\begin{lemma}\label{le:turn1}  Let $\phi_{t \in [0,1]}$ be an essential
loop in $\Ham(S^2)$, where $S^2$ has area $A$.
Then $\Ll(\phi_{t \in [0,1]}) \ge A$.
\end{lemma}

\proof{}
Let $F:S^2 \to \R$ be the composite of the height function on $S^2$
with a linear map so that $\|F\| = A$.
It is easy to check that its flow is exactly a full rotation
about the poles. Since its orbits all have period exactly $1$,
Proposition~\ref{prop:ball4} shows that, for any $\eps > 0$,
one can embed a ball of capacity $A-\eps$ on both sides of
$gr(F)$.  (In fact, this is another description of the full
filling of $S^2\times D^2$ by $2$ balls which was given in
\cite{TRAY}.)  Thus, by Theorem~\ref{thm:monotone},
$$
\Ll(\phi_{t \in [0,1]}) \geq A.\eqno\Box
$$

\begin{remark}\label{rmk:cap}\rm
Proposition~\ref{prop:ball4} shows that in dimension $2$ the
condition ``$c_G(H) = \|H\|$'' is no weaker than the condition
``$H$ has no non-constant periodic orbit in time $< 1$''.
Actually,
it is strictly stronger, even on $\R^2$.  Consider any positive
(autonomous) bump function $H$ with support of
area less than half the area of $M$, but with ${\rm area}(H^{-1}(\max H))
> \max H$.  Then
 $c_G(H) = \|H\|$ because one can then embed a ball of capacity $\|H\|$ 
on both sides
of the graph of $H$, and our argument shows that the flow $\phi_{t \in [0,1]}$ of
$H$ is length-minimizing. But  $H$ does have non-constant periodic orbits 
in time less than $1$
as soon as the slope of the bump near the boundary of its support is high
enough. 
\end{remark}

\subsection{Hofer's diameter}

We briefly discuss questions concerning the diameter of 
$\Ham^c(M)$ under the Hofer norm $\|\cdot \|$.  First observe that if $(M,\om)$ is
exact and of finite volume, it follows easily from the existence of the Calabi
homomorphism
$$
\Cal: \Ham^c(M) \to \R
$$ 
that this diameter is infinite. 

   Now the argument in Remark~\ref{rmk:cap}
 shows that the diameter of $\Ham^c(M)$, for
a geometrically bounded semi-monotone manifold $M$ with $r_1 = \infty$ 
containing 
two disjoint embedded copies of $B^{2n}(c)$, is at least equal to $c$.
In some cases it is enough to know that a covering of $M$ contains such
balls.  The next lemma shows that the diameter of $\Ham(T^2)$ is infinite.
For simplicity, the argument uses the energy-capacity inequality, but it
could equally well be phrased in terms of embedded balls.

\begin{lemma}  There is a function $H$ on $T^2$ whose flow $\{\phi_t^H\}_{t\ge
0}$ is such that
$$
\Ll(\phi_t^H) = t \, {\rm TotVar}(H), \;\;{\rm for} \; {\rm all} \; t \ge 0.
$$
 \end{lemma}

\proof{}  Consider a Hamiltonian of the form $H(x,y) = f(x)$,
where 
\begin{description}
\item[(i)] $f(0) = f'(0) = 0$, and

\item[(ii)]  $f$ increases
over $0\le x\le 1/2$ to $f(1/2) = 1$ and then decreases to $0$ at $x= 1=0$. 
\end{description}
If 
$\phi_t, t\ge 0$, is its flow,  let $\Tilde\phi_t$  be the unique lift of $\phi_t$
to $\R^2$ which starts at $\1$.  Then
$$
\Tilde\phi_t(x,y) = (x, y+ tf'(x)). 
$$
Hence in time $T$ this isotopy
disjoins a region $U_+$ in the strip  $0< x < 1/2$, which is diffeomorphic to a
disc and has area almost equal to $T(f(\frac 12) - f(0)) = T$ and a similar region
$U_-$ in the strip $1/2 < x < 1$.   

Let $K_t$ be any
Hamiltonian on $T^2$ with time-$1$ map $\phi_T$.  Then its lift $\Tilde K_t$ to
$\R^2$ has time $1$-map $\Tilde\psi = \Tilde\phi_T + c$ for some $c\in \Z^2$.
Since $\Tilde\phi_T$ moves $U_+$ upwards and $U_-$ downwards, $\Tilde\psi$ must
disjoin at least one of these.  But  the energy-capacity inequality
of~\cite{HOF,LALMCD} states that any map $\psi$ which disjoins a disc of area $c$
must have norm $\|\psi\|$ at least $c$.  Thus
$$
\Ll(\{\Tilde K_t\}) = \Ll(\{K_t\}) \ge T,
$$
as required. \QED

\begin{cor}  The norm on $\Ham (T^2)$ is unbounded.
\end{cor}

A similar argument clearly applies to any compact Riemann surface $\Sigma$  
except
the sphere.  (Note that the universal cover of such $\Sigma$ is still
symplectomorphic to the plane $\R^2$: its negative curvature is not reflected in
its symplectic structure.)  The argument also  works for manifolds of the form
$\Sigma\times M$, where $M$ is any compact manifold whose universal cover has
infinite Gromov capacity, since in this case again one can lift $\phi_t$ to an isotopy which
disjoins an arbitrarily large ball.  

 In fact, the methods of~\cite{LALMCD} show that one can prove a similar result  for 
$\Si\times M$ for any
compact $M$ provided one can show that there is a 
symplectic embedding 
$$ 
\iota: B^4(r)\times M \to B^2(R)\times \R^2\times M.
$$
only if $r \le R$.
This generalized non-squeezing theorem can be proved by adapting the methods
of~\cite{LALMCD}.

\subsection{Linear rigidity}

Siburg showed in~\cite{S} that when $M = \R^{2n}$
any path $\phi_{t\in [0,1]}$ such that
$\phi_t$ has no non-constant closed trajectories for any $t\in [0,1]$ is
absolutely length-minimizing.  (Such non-constant trajectories correspond
precisely to fixed points of $\phi_t$ which are non-trivial in the sense that
 $\phi_s(x) \ne x$ for some $s\in [0,t]$.)
   Combining this with our necessary condition for a stable geodesic in
Hofer's metric (see \cite{LALM1})
we obtain:

\begin{theorem}[Linear rigidity] \label{thm:linear-rig} 
Let $H_{t \in [0,1]}$ be a quasi-autonomous
Hamiltonian  on $R^{2n}$, with say only one fixed minimum 
$p \in \R^{2n}$. If the linearised flow $d\phi_t(p)$ has a non-constant closed
trajectory in time less than $1$, so does the flow $\phi_t$ itself. 
The same conclusion still
holds if  $H_{t \in [0,1]}$ has a finite number of fixed extrema 
provided that 
the linearised flows at all fixed minima or at all fixed maxima have  
a non-constant closed trajectory in time less than $1$.
\end{theorem}

\proof{} If the conclusion does not hold, 
then the path $\phi_{t \in [0,1- \eps]}$
is length-minimizing by Siburg's criterion, for any choice of small $\eps > 0$.
It must then be a stable geodesic,
which implies the non-existence of non-constant closed trajectories of the
linearised flow in time $< 1 - \eps$ by the necessary condition for stability of
\cite{LALM1}. Since $\eps$ is arbitrarily small, this contradicts the hypothesis.
\QED   

     Note that we need no non-degeneracy condition of $H_t$ at $p$. Actually, this
theorem is non-trivial even in the autonomous case.

   Thus the linearised flow controls the flow itself, or at least the behaviour
of its closed trajectories. The above proof makes use of Siburg's
criterion only to show that the given path is length-minimizing in its homotopy 
class. Thus a similar result would hold on any weakly exact manifold for which
one can prove Proposition~\ref{prop:HZ-Hamiltonians} for quasi-autonomous
Hamiltonians.

   We do not know whether the closed trajectory given by the theorem can be
found in any arbitrarily small neighbourhood of $p$.

   Of course, this linear rigidity is a symplectic phenomenon: on
$\R^2$, there are non-symplectic autonomous flows having no non-constant
closed orbit and whose linearization at $0$ is a pure rotation.

\subsection{$C^1$-flatness on  manifolds without short loops}

We will now calculate the Hofer norm on a neighbourhood  $\Uu$ of the identity
in $\Ham^c(M)$ assuming that $M$ satisfies the conditions of the first part of
Theorem~\ref{thm:lengmin} ($M$ is semi-monotone and has no short loops).
  As in Lemma~\ref{le:split}, we  take $\Uu$ to be  a
star-shaped neighbourhood of $\1$ in $\Ham^c(M)$ consisting of  Hamiltonian
diffeomorphisms $\psi$ whose graphs lie close enough to the diagonal $diag$
in $(M\times M, -\om\oplus \om)$ to correspond to an exact $1$-form
$\rho(\psi)$ in $(T^*M, -d\,\la_\can)$.\footnote
{
These sign conventions are the same as in~\cite{DUSAL2}.  More details
 of the arguments given here may be found in Chapter~9 of that book.
}
Then there is a function  $F_\psi$ on $M$ which is unique up to a constant and
such that $\rho(\psi) = dF_\psi$, and we claim:

\begin{proposition}\label{prop:BP}  Consider the family of $1$-forms $t
dF_\psi,\, 0\le t \le 1$, and let $\psi_t$  be the corresponding isotopy in
$\Ham^c(M)$.  Then, if $F_\psi$ is sufficiently $C^2$-small, 
$$
\Ll(\psi_{t\in [0,1]}) = \sup_x F_\psi - \inf_x F_\psi
$$
and each extremum of $F$ is a fixed extremum of $H_t$. Here $H_t$ is the
Hamiltonian generating $\psi_t$.
\end{proposition}

\proof{}  Observe that  the fixed points of $\psi$ correspond precisely to the
critical points of $F = F_\psi$ and hence are fixed by $\psi_t$ for all $t$.
Therefore, if $F$ assumes its maximum and minimum values at $P,p$, these
points are fixed points of the isotopy $\psi_t$.  Let $H_t$ be the Hamiltonian
which generates $\psi_t$, and $\be(s), 0\le s \le 1$, be any path in $M$ from
$\be(0) = p$ to $\be(1) = P$.  Then we claim that both
$\Ll(\psi_t) $ and $ \sup_x F_\psi - \inf_x F_\psi = F(P) - F(p)$ are equal to
the area  $A$ swept out by the arc $\be$ under the isotopy $\psi_t$.  
More
formally, this area can be written as
$$
\int_{[0,1]\times [0,1]} \Psi^*(\om),
$$
where $\Psi(s,t) = \psi_t(\be(s))$.

We will prove this in a slightly different form.  Given any two critical points
$q_1, q_2$ of $F$, let $\be$ be any arc in $M$ from $q_1$ to $q_2$ and write
$A(q_1,q_2)$ for the area swept out by $\be$ under the isotopy $\psi_t$.
Then we will show that
$$
F(P) - F(p) = \sup_{q_1,q_2} A(q_1,q_2) =\Ll(\psi_t).
$$

  To prove the statement about $F$ observe that the path
$(\be(s),\psi_t(\be(s))$ on the graph of $\psi_t$ corresponds to a path in
$T^*M$  from $(q_1,0)$ to $(q_2,0)$ which lies on the section of $T^* M$
determined by the $1$-form $tdF_\psi$.  Thus as $t$ varies in $[0,1]$ they
form a homotopy of paths (with fixed endpoints) in $T^*M$ from the path
$\be$ in the zero-section to the path $\psi(\be)$ in graph $\,(dF)$. If we
write $Y'$ for the corresponding $2$-cycle in $T^*M$, we find that, because
$\la_\can = dF$ on  graph $\,(dF)$,
 $$
\int_{Y'} -d\,\la_\can  = F(q_2) - F(q_1).
$$
Now, look at the corresponding set $Y'' $ traced out by the paths
$(\be(s),\psi_t(\be(s))$ in $(M\times
M,-\om\oplus\om)$.
Then, because  $Y''$  lies over a $1$-dimensional subset of
the first factor of $M$, namely $\be$, we find,
 $$
\int_{Y'} -d\,\la_\can  =\int_{Y''} -\om\oplus\om =  \int_{[0,1]\times
[0,1]}\Psi^*(\om) = A(q_1,q_2).
$$
Since $F(P) - F(p) = \sup_{q_1,q_2}(F(q_2) - F(q_1))$, the result follows.
(Note that this calculation shows that $A(q_1,q_2)$ is independent of the
choice of $\be$.)

 To prove the statement for $\Ll(\psi_{t\in [0,1]})$, consider
the graph of $H_t$: 
$$
\Ga_H = \{(x,s,t)\in M\times \R\times [0,1] \mid  s = H_t(x)\}.
$$
Let $Y\subset \Ga_H$ be the union of all the characteristic lines in $\Ga_H$
starting at the points of $\be$:
$$
Y = \{ (\psi_t\be(s),  H_t(\psi_t\be(s)), t)\in \Ga_H:s,t\in [0,1]\}.
$$
Then $\om \oplus \si$ vanishes on $Y$ because $Y$ is a union of null lines.
Hence
$$
 \int_Y \om =  - \int_Y \si.
$$
But the first integral here is just $A(q_1,q_2)$, while the second is minus the
area of the projection of $Y$ onto $\R^2$, i.e. minus the area
$$
\area_{t\in [0,1]}\{q_1,q_2\}
$$
 enclosed
by the curves $H_t(q_2)_{t\in [0,1]}$ and $H_t(q_1)_{t\in[0,1]}$.  Thus we
have shown that for any pair of critical points of $F$
$$
F(q_2) - F(q_1)= A(q_1,q_2) = -\area_{t\in [0,1]}\{q_1,q_2\}.
$$
Note that the corresponding inequality holds over any time interval
$t\in[a,b]\subset[0,1]$.   

It is clear from the definition of
$H_t$ that for each $t$ this function has the same critical points as $F$. 
 (These are
the infinitesimally fixed points of the isotopy $\psi_t$.) But the fact that
the equality $(b-a)(F(q_2) - F(q_1))=-\area_{t\in [a,b]}\{q_1,q_2\}$ holds
for all $0 \leq a \leq b \leq 1$ implies easily that any minimum $q_t$ of
some $H_t$ must be a minimum of $F$ too, and therefore a fixed minimum
of the family $H_{t \in [0,1]}$. This is because the function $F$ is
independent of time and the left hand side of the above equation is maximized
at the extema of $F$. Since the same applies to maxima, this 
proves the last assertion of the Proposition. The first assertion follows
at once:
$$
\Ll(\phi_t) = \area_{t\in [0,1]}\{q, Q\} =  
-A(P,p) = F(P) - F(p),
$$
where $q$ and $Q$ are fixed minimum and maximum of $H_t$ and $P,p$ are
the corresponding maximum and minimum  of $F$.  \QED

\begin{cor}  If $M$ is semi-monotone and has no short loops and if
the neighborhood $\Uu$ is sufficiently $C^1$-small, then
$$
\|\psi\| = F(P) - F(p).
$$
\end{cor}
\proof{}  It is easy to check that if $F$ is $C^2$-small so is the generating
Hamiltonian $H_t$ of the isotopy $\psi_t$.   By Lemma~\ref{le:C2}, if
$H_t$ is sufficiently $C^2$-small then $\Ll(H_t) = c_G(H)$.  Hence, by
Theorem~\ref{thm:lengmin}, the path $H_t$ measures the length of $\psi$.
Thus
$$
\|\psi\| = \Ll(H_t) = F(P) - F(p).\eqno\Box
$$

\begin{remark}\rm  One of the key points in the proof of
Proposition~\ref{prop:BP} is to understand the relation between the generating
function $F_\psi$ and the Hamiltonians  $H_t$.  This is discussed further
in~\cite{DUSAL2}.  We have argued geometrically and somewhat indirectly here, but
one could construct a proof along the lines of that in Bialy--Polterovich, using
an appropriate generalization of the Hamilton-Jacobi equation.
\end{remark}

\end{document}